\documentclass[aap]{imsart}

\RequirePackage{amsthm,amsmath,amsfonts,amssymb}
\RequirePackage[numbers]{natbib}
\RequirePackage[colorlinks,citecolor=blue,urlcolor=blue]{hyperref}
\RequirePackage{graphicx}

\usepackage{mathtools, bm, bbold}
\usepackage[OT1]{fontenc}
\usepackage{parskip}
\usepackage{enumitem}
\usepackage[ruled]{algorithm2e}
\usepackage{tabu,booktabs}
\usepackage{subcaption,caption}

\startlocaldefs
\numberwithin{equation}{section}
\theoremstyle{plain}
\newtheorem{theorem}{Theorem}[section]

\theoremstyle{remark}
\newtheorem{remark}[theorem]{Remark}
\newtheorem{definition}[theorem]{Definition}

\newtheorem{property}[theorem]{Property}
\newtheorem{dataset}[theorem]{Dataset}

\newif\ifinthesis
\inthesisfalse
\newcommand{\figuretb}{\begin{figure}[tb]}
\newcommand{\tabletb}{\begin{table}[tb]}
\newcommand{\algohtb}{\begin{algorithm}[!htb]}
\newcommand{\algotb}{\begin{algorithm}[tb]}
\newcommand{\figurehtb}{\begin{figure}[!htb]}

\renewcommand{\phi}{\ensuremath{\varphi}}
\renewcommand{\epsilon}{\ensuremath{\varepsilon}}
\DeclareMathOperator*{\argmin}{arg\,min}
\DeclareMathOperator*{\BigCup}{\bigcup}

\newcommand{\pluseq}{\mathrel{{+}{=}}}
\newcommand{\moinseq}{\mathrel{{-}{=}}}

\newcommand{\xdata}{\mathcal X}
\newcommand{\xdatai}{\bm X_i}
\newcommand{\xdataN}{\bm X_1,...,\bm X_N}

\endlocaldefs

\begin{document}

\begin{frontmatter}
\title{Estimation of high dimensional Gamma convolutions through random projections}
\runtitle{Estimation of high dimensional Gamma convolutions}

\begin{aug}
\author[A]{\fnms{Oskar} \snm{Laverny}\ead[label=e1]{oskar.laverny@univ-lyon1.fr}},
\address[A]{Institut Camille Jordan UMR 5208, Université Claude Bernard Lyon 1, France,\\ and SCOR SE, \printead{e1}}
\end{aug}

\begin{abstract}
    Multivariate generalized Gamma convolutions are distributions defined by a convolutional semi-parametric structure. 
Their flexible dependence structures, 
  the marginal possibilities and their useful convolutional expression 
  make them appealing to the practitioner. 
However,
  fitting such distributions when the dimension gets high is a challenge.
We propose stochastic estimation procedures based on the approximation 
  of a Laguerre integrated square error via (shifted) cumulants approximation,
  evaluated on random projections of the dataset.
Through the analysis of our loss via tools from Grassmannian cubatures, sparse optimization on measures and Wasserstein gradient flows, 
  we show the convergence of the stochastic gradient descent to a proper estimator
  of the high dimensional distribution.
We propose several examples on both low and high-dimensional settings.
\end{abstract}



\end{frontmatter}

{\graphicspath{{img/}} 
\section{Introduction}
\ifinthesis\addtocontents{toc}{\protect\vskip-11pt}\fi


The class of univariate generalized Gamma convolutions was introduced by Thorin~\cite{thorin1977,thorin1977a} in 1977 to prove the infinite divisibility of the log-Normal and Pareto distributions. 
Bondesson~\cite{bondesson1992} formalized in the early 1990 a large portion of the theory, his book is still a great resource nowadays. 
Thorin's class is the class of weak limits of convolutions of independent Gamma distributions. 
It is obviously closed by convolution, but also by independent product of distributions, exponentiation, and many other interesting applications, as Bondesson showed in~\cite{bondesson2015} recently. 
It contains many distributions that are useful to the practitioner, such as the log-Normal and Pareto as historically important examples, but also some Weibulls, $\alpha$-stables distributions, and all product, sum and exponentiation of those. 

The multivariate Thorin class was introduced by Bondesson in~\cite{bondesson2009}, extending an old idea from Cherian~\cite{cherian1941}. Sadly, the theoretical characteristics of the class were problematic and the extension of the univariate theory was left almost untouched. Noteworthy is the work from Perez \& Abreu~\cite{perez-abreu2012,perez-abreu2014} that generalized the concept to other cones than $\mathbb R_+^d$. 

From a statistical point of view, this class is appealing for several reasons. First, the marginals belong to the well-studied univariate Thorin class, and the multivariate case allows for a wide range of dependence structures between them. Second, the marginal distributions can be heavy-tailed, and do not need to have moments. Last, the whole multivariate distribution has a parametrization that is interpretable and simple, by definition infinite divisible, and which simplifies deconvolutions computations and the analysis of multivariate extremes. 

While the univariate Thorin class benefits from a large description of its content and properties, e.g., in~\cite{furman2017,miles2019}, or in~\cite{roynette2009,james2008,bondesson2015,bondesson2018}, the statistical literature and especially regarding the problem of the estimation of such parametric models is a lot sparser. One potential reason is that the deconvolution of an empirical distribution into several Gamma distributions is a hard inverse problem, and the deconvolution literature already has troubles with removing a small convolutional noise from a strong signal.

Miles, Furman and Kuznetsov~\cite{furman2017,miles2019} provide a way to project a known density onto the univariate class, but the procedure is somewhat complicated and requires very precise integration of shifted moments. Up until~\cite{laverny2021a} recently, there were no estimation procedure for distributions in this class, neither univariate nor multivariate. This estimation procedure uses a projection in a Laguerre basis~\cite{dussap2021,mabon2017,comte2015}, together with a bijection between parameters of the distribution and its coefficients in the basis, to produce an efficient loss and a fast numerical procedure to fit multivariate generalized Gamma convolutions. However, the procedure is not adapted to high dimensional cases as it suffers badly from the dimensionality curse.   

We derive here a loss and a corresponding estimator that uses random projections, exploiting the convolutional properties of the class, to solve the dimensionality issues. With this loss, we are able to fit correctly multivariate Gamma convolutions on datasets in dimension up to $2000$ in a few minutes on a standard laptop, as our examples show.

The paper is organized as follows: 
\begin{itemize}
  \item Section~\ref{ggchd:sec:previous_work} and Section~\ref{ggchd:sec:random_projs} give definitions of the generalized Gamma convolutions, fix some notations, describe the approach from~\cite{laverny2021a} which we leverage, and discuss briefly Grassmannian cubatures and random projections of moments problems.
  \item Section~\ref{ggchd:sec:convex_approx} proposes a new loss to fit generalized Gamma convolutions, which we show to be consistent under well-behavior assumptions, and Section~\ref{ggchd:sec:gradient_flows} establishes the convergence of the gradient flow of our loss to a global minimizer in the space of Thorin measures.
  \item Section~\ref{ggchd:sec:investigation} contains investigation of the proposed estimator and a few numerical examples in various dimensions, and Section~\ref{ggchd:sec:conclusion} concludes.
\end{itemize}

\section{Multivariate Thorin measures}\label{ggchd:sec:previous_work}
\ifinthesis\addtocontents{toc}{\protect\vskip-11pt}\fi

We start by a quick remark on notations. 
\begin{remark}[Notations] 
  We use bold letters such as $\bm a$ to designate generally (finite dimensional) indexable objects,
    e.g., vectors or matrices,
    and corresponding indexed and possibly unbolded letters,
    such as $a_i$, for an index set $I$ s.t. $\lvert I \rvert < \infty$,
    designate values in these objects.
  Inequalities, additions, products, and fractions between vectors of same shape are intended componentwise, and standard broadcasting works with scalars so that $$2\frac{\bm x}{\bm y} \le \bm z + \bm t \bm u \iff \forall i \in I,\;2\frac{x_i}{y_i} \le z_i + t_i u_i.$$ On the other hand, powers, norms and scalar products are reducing operations so that:
   $$\bm x^{\bm y} = \prod\limits_{i\in I} x_i^{y_i}\text{, } 
      \langle\bm x,\bm y\rangle = \sum\limits_{i\in I} x_iy_i  \text{ and } 
      \lVert \bm x\rVert_2^2 = \langle \bm x, \bm x\rangle$$
  where $\langle \cdot ,\cdot \rangle$ denotes the Euclidean scalar product on $\mathbb R^{\lvert I \rvert}$. When $\bm A$ is a matrix with column indices in a set $J$ and row indices in $I$, $\bm A\bm x = \left(\sum_{j \in J} A_{i,j}x_j\right)_{i \in I}$ denotes the classical matrix-vector product.
\end{remark}

For the exposition of the whole paper, we consider $\bm X = \left(X_1,...,X_d\right)$ to be a random vector in $\mathbb R_+^d$. 
Our goal is to assess the distribution of $\bm X$.
Typically, 
  for an insurer, 
  $X_1,...,X_d$ might be losses from several lines of business.
Recall that such a random vector can be characterized by its cumulant generating function 
\begin{align*}
  K\colon \mathbb C^d   &\to \mathbb C\\
      \bm t       & \mapsto \ln\mathbb E\left(e^{\langle \bm t, \bm X \rangle}\right),
\end{align*}
  where it is defined 
  (it is always defined at least around $\bm t = -\bm 1$, and more generally in the negative half of the space since $\bm X$ is positive). Bondesson~\cite{bondesson1992,bondesson2009} defined the class of multivariate Gamma convolutions, letting $\mathcal M_+(\mathbb R_+^d)$ be the set of positive measures with support $\mathbb R_+^d$, as follows: 

\ifinthesis \restatedefGd \else 
\begin{definition}[$\mathcal G_{d}$]\label{ggchd:def:Gd}
A random vector $\bm X$ is supposed to follow a $d$-dimensional generalized Gamma convolution
  with Thorin measure $\nu \in \mathcal M_+(\mathbb R_{+}^d)$,
  hereafter denoted $\bm X \sim \mathcal G_{d}(\nu)$, 
  if its cumulant generating function writes: 
  $$K(\bm t) = - \int_{\mathbb R_+^d} \ln\left(1 - \langle \bm t,\bm s \rangle\right) \nu(d\bm s).$$

When $\nu$ is finitely atomic, i.e., denoting $\delta_{x}$ the Dirac measure in $x$, when there exists $n \in \mathbb N$, $\bm\alpha \in \mathbb R_+^n$ and $\bm s \in \mathbb R_+^{n\times d}$ such that $\nu = \sum\limits_{i=1}^n \alpha_i \delta_{\bm s_i}$, $\bm X$ is said to be a $d$-dimensional Gamma convolution with shapes $\bm\alpha$ and scales $\bm s$.
\end{definition}
\fi

Not every $\nu \in \mathcal M_+(\mathbb R_{+}^d)$ is the Thorin measure of a multivariate generalized Gamma convolution, as the resulting function $K$ from Definition~\ref{ggchd:def:Gd} might not be a cumulant generating function. 
Conditions on $\nu$ making $K$ a proper cumulant generating function can be found in~\cite{bondesson1992} for $d=1$ and in~\cite{perez-abreu2012,perez-abreu2014} for $d > 1$. Of course, this is always the case for finitely atomic $\nu$. 

Note that in the finitely atomic case, say $\bm X \sim \mathcal G_{d}(\sum_{i=1}^n \alpha_i \delta_{\bm s_i})$, there exists independent Gamma random variables $G_{i} \sim \mathcal G_{1}(\alpha_i\delta_{1})$, all having unit scale, such that:

\begin{equation*}
  \begin{pmatrix}
      X_{1} \\
      ... \\
      X_{d}
  \end{pmatrix}
    =
  \begin{pmatrix}
      s_{1,1} & ... & ... & s_{1,n} \\
      ... & ... & ... & ... \\
      s_{d,1} & ... & ... & s_{d,n}
  \end{pmatrix}
    \cdot
  \begin{pmatrix}
      G_{1} \\
      ... \\
      G_{n}
  \end{pmatrix},
\end{equation*}

where typically many $s_{i,j}$ are zeros. To fix the parametrization in the readers mind, remember that the above vector $\bm X$ has mean $\mathbb E\left(\bm X\right) =\bm s\bm\alpha$ and variance-covariance matrix $\mathbb V\left(\bm X\right) = \bm s\mathrm{diag}(\bm\alpha)\bm s$.

This convolutional structures of distributions in $\mathcal G_d$ is appealing for a few different reasons. First, it allows computing directly the distribution of a linear aggregation of the random vector, even under dependency. Second, it allows for parametric (and thus fast) divisibility\footnote{A random vector $\bm X$ is said to be $n$-divisible if we can find independent and identically distributed random vectors $\bm X_1,...,\bm X_n$ independently distributed such that $\bm X = \bm X_1 +...+\bm X_n$. It is said to be infinitely divisible if this can be done for any $n \in \mathbb N$.}, which can be of importance for certain applications. Third, the class of achievable marginals is quite large and contains many interesting distributions on the positive real line with light or heavy tails, finite moments or not, etc. Fourth, (non-generalized) multivariate Gamma convolutions are really easy to sample and work with. If they are estimated on real datasets, the sparsity structure of the matrix $\bm s$ can be valuable information to the practitioner, which could lead to insight about the dependence structure underlying the dataset.  

The wide cases ($n \gg d$) and tall cases ($d \gg n$) are both of interests: they can be seen as (resp) under determined and over determined non-negative blind source separations problems (see~\cite{comon2010}), where the sources are supposed to be independent and Gamma distributed. 

In-depth study of distributions in $\mathcal G_1$ (i.e., achievable marginals in $\mathcal G_d$) is done in~\cite{bondesson1992}. 
Mataï~\cite{mathai1982} and later Moschopoulos~\cite{moschopoulos1985} 
  provide convergent series to estimate the densities of $\mathcal G_{1}$ distributions with finite Thorin measure.  
In~\cite{laverny2021a},
  a more systematic approach is given to the multivariate density problem through Laguerre expansions, generalizing Moschopoulos's series to a multivariate case. 

\begin{definition}[Laguerre basis~\cite{laverny2021a,dussap2021,comte2015,mabon2017}]
The Laguerre functions $\left(\phi_{\bm k}\right)_{\bm k \in \mathbb N^d}$,
  defined $\forall\bm x \in \mathbb R_+^d$ by 
  $$\phi_{\bm k}(\bm x) = \prod_{i=1}^d \phi_{k_i}(x_i) \text{ for }\bm k \in \mathbb N^d,
  \text{ where } 
  \phi_k(x) = \sqrt{2} \sum_{j=0}^{k} \binom{k}{j} \frac{(-2x)^{j}}{j!} e^{-x} \text{ for }k \in \mathbb N,
  $$
  form an orthonormal basis of the space $L^2(\mathbb R_{+}^d)$ of square integrable functions on $\mathbb R_+^d$.
\end{definition}

\begin{definition}[Laguerre coefficients $a_{\bm k}$~\cite{dussap2021}] When the density $f$ of the random vector $\bm X$ is square integrable, it can be expanded in the Laguerre basis, and we denote $a_{\bm k}$ its coefficients, i.e., $$f(\bm x) = \sum\limits_{\bm k \in \mathbb N^d} a_{\bm k} \phi_{\bm k}(\bm x) \text{ and therefore } a_{\bm k} = \mathbb E\left(\phi_{\bm k}(\bm X)\right)$$
\end{definition}

\begin{definition}[Shifted moments $\mu_{\bm k}$] We denote $\mu_{\bm k} = \mathbb E\left(\bm X^{\bm k} e^{\langle -\bm 1, \bm X \rangle}\right)$ the $\bm k^{\text{th}}$ shifted moment of $\bm X$, defined as the $\bm k^{\text{th}}$ Taylor coefficient of the moment generating function $M$ of the random vector $\bm X$ expanded around $-\bm 1$: 
  $$M(\bm t) = \mathbb E\left(e^{\langle \bm t, \bm X\rangle}\right) 
  = \sum\limits_{\bm k \in \mathbb N^d} \frac{\mu_{\bm k}}{\bm k!} \left(\bm t - \bm 1\right)^{\bm k},$$

  (where $\bm k ! = \prod_{k \in \bm k} k!$). 
\end{definition}

Due to the linearity of the expectation operator and the definition of $\phi_{\bm k}$, we have the following linear link between Laguerre coefficients and shifted moments:
  \begin{equation}\label{ggchd:eq:relation_a_mu}
    a_{\bm k} = \sqrt{2}^d \sum_{\bm j\le \bm k} \binom{\bm k}{\bm j} \frac{(-2)^{\lvert\bm j\rvert}}{\bm j !} \mu_{\bm j}.
  \end{equation}

Let us now pack up these shifted moments and Laguerre coefficients and reframe this relation as a matrix-vector expression. 
%
%
For all the exposition, let $I$ be a (non-empty) increasing index set, that is a finite set of indices $I \subset \mathbb N^d$ such that if $\bm i \in I$, 
then each $\bm j \le \bm i$ (componentwise) is also in $I$. Pack up Laguerre coefficients and shifted moments as: $$\bm a := \left(a_{\bm k}\right)_{\bm k \in I},\text{ and }\bm \mu := \left(\mu_{\bm k}\right)_{\bm k \in I}.$$

Remark that the relationship between $\bm a$ and $\bm \mu$ induced by Equation~\eqref{ggchd:eq:relation_a_mu} is linear,
  bijective,
  and does not depend on the distribution of $\bm X$.
We encode this relationship through a (lower-triangular, invertible) matrix
  $\bm A \in \mathbb R^{\lvert I\rvert \times \lvert I\rvert}$ 
  such that $$\bm a = \bm A\bm \mu \text{ and }\bm \mu = \bm A^{-1} \bm a.$$ 

We now consider the derivatives of the cumulant generating function $K = \ln \circ M$. For $\bm p \in \mathbb N^d$, we denote the $\bm p^{\text{th}}$ multivariate derivatives of a function $f$ as $f^{(\bm p)}(\bm t) = \frac{\partial^{\lvert \bm p \rvert}}{\partial^{p_1} t_1 \cdot ... \cdot \partial^{p_d} t_d} f(\bm t)$. Then we define the Thorin moments as follows:

\begin{definition}[Thorin moments $\tau_{\bm k}$]\label{ggchd:def:tau} Let $\tau_{\bm k} = \frac{K^{(\bm k)}(-\bm 1)}{(\lvert\bm k\rvert -1)!}$ be the $\bm k^{\text{th}}$ Thorin moment of the random vector, i.e., the constant such that the expansion of the cumulant generating function $K$ around $-\bm 1$ writes:
  $$K(\bm t) = \ln M(\bm t) = \sum\limits_{\bm k \in \mathbb N^d} \frac{\tau_{\bm k}(\lvert\bm k\rvert -1)!}{\bm k!} \left(\bm t - \bm 1\right)^{\bm k}.$$ 
\end{definition}

Reasons that push us to define $\tau_{\bm k}$ as it is and not as directly the Taylor coefficients of $K$ at $-\bm 1$ will appear shortly in Property~\ref{ggchd:prop:thorin_moments}. Once again, we pack them as: $$\bm \tau := \left(\tau_{\bm k}\right)_{\bm k \in I}.$$

Since $M(\bm t) = e^{K(\bm t)}$, there is also a bijection between $\bm \tau$ and $\bm \mu$, which we encode through a function $\bm B$ such that 
$$\bm \mu = \bm B(\bm \tau) \text{ and } \bm \tau = \bm B^{-1}(\bm \mu).$$

The function $\bm B$ is a complicated (but continuous, bijective and smooth) function, mostly relying on multivariate Bell's polynomials. See Appendix~\ref{ggchd:apx:nabla_AB} for details on the function $\bm B$ and for computationally fast recursive formulation of the function $\bm B$. 

\begin{remark}[Expansion point] We expanded $M$ and $K$ around $-\bm 1$ and not around $\bm 0$. Indeed, expansion around $\bm 0$, although more practical, requires hypothesis on the random vector that we are not willing to make (mainly moments hypothesis). Here, we only supposed that $\bm X$ is positive. Of course, any other expansion point $\bm t \in \mathbb C^d$ such that each $\mathrm{Re}(t_i) <0$ will do, both functions are analytic in the negative half of the space and the Laguerre basis can be rescaled accordingly (see~\cite{mabon2017}). To simplify the exposition, we chose $-\bm 1$. 
\end{remark}

Let us now discuss the particular case of $\mathcal G_d$ distributions. 

\begin{definition}[Thorin moments of a Gamma convolution]\label{ggchd:def:tau_of_nu} We denote the first Thorin moments (from Definition~\ref{ggchd:def:tau}) of the $\mathcal G_{d}(\nu)$ distribution by: $$\bm \tau(\nu) = \left(\tau_{\bm k}(\nu)\right)_{\bm k \in I}.$$
\end{definition}

As we already noted, the definition of the Thorin moments includes a renormalization of standard Taylor coefficients of $K$. This is done so that Property~\ref{ggchd:prop:thorin_moments} holds.
\begin{property}\label{ggchd:prop:thorin_moments} Denoting $\delta_x$ the Dirac measure at $x$, it holds: 
  \begin{enumerate}[label=(\roman*)]
    \item $\tau_{\bm 0}(\delta_{\bm s}) = -\ln\left(1 + \lvert \bm s\rvert\right)$
    \item $\tau_{\bm k}(\delta_{\bm s}) = \bm s^{\bm k}\left(1 + \lvert \bm s\rvert\right)^{-\lvert \bm k \rvert}$ for $\bm k \neq 0$, $\bm k \in \mathbb N^d$.  
    \item $\bm \tau(\nu) = \int \bm \tau(\delta_{\bm s}) \nu(d\bm s)$, where the equality and integration are intended componentwise.  
  \end{enumerate}
  \begin{proof} 
    The first point is simply the definition of the cumulant generating function of the $d$-variate Gamma distribution $\mathcal G_{d}(\delta_{\bm s})$, taken at $\bm t=-\bm 1$.
    The third point comes directly from Definition~\ref{ggchd:def:tau_of_nu}, and from the fact that cumulants are linear on convolutions of independent random vectors.
    To get the second point, derive the cumulant generating function of the $\mathcal G_{d}(\delta_{\bm s})$ distribution as follows:
  \begin{align*}
    \tau_{\bm k}(\delta_{\bm s}) 
    &= \frac{1}{(\lvert\bm k\rvert -1)!}\left. \frac{\partial^{\bm k}}{\partial^{\bm k} \bm t} K(\bm t) \right\rvert_{\bm t = -\bm 1}\\
    &= \frac{1}{(\lvert\bm k\rvert -1)!}\left. \frac{\partial^{\bm k}}{\partial^{\bm k} \bm t} -\ln\left(1 - \langle \bm t,\bm s \rangle\right) \right\rvert_{\bm t = -\bm 1}\\
    &= \frac{1}{(\lvert\bm k\rvert -1)!}\left. \bm s^{\bm k} (\lvert \bm k\rvert -1)! \left(1 - \langle \bm t, \bm s\rangle\right)^{-\lvert \bm k \rvert} \right\rvert_{\bm t = -\bm 1}\\
    &= \left(\frac{\bm s}{1 + \lvert \bm s\rvert}\right)^{\bm k}.
  \end{align*}
  \end{proof}
\end{property}

We note that the Thorin moments are not exactly the moments of $\nu$, but rather generalized moments of $\nu$. However, using the change of variables $\bm x =  \frac{\bm s}{1 + \lvert\bm s\rvert}$ in the integral from Property~\ref{ggchd:prop:thorin_moments}, point $(iii)$, allows to express them as the moments of another measure $\xi$. See~\cite{james2008} for a large discussion of the $\bm x, \xi$ parametrization versus the $\bm s, \nu$ parametrization in the univariate case, 
  whose motivations also apply to the multivariate case. 
This is a generalization of the rate-or-scale dilemma for Gamma distributions parametrizations. 
Others parametrizations are possible, 
  see~\cite{perez-abreu2014} for one that allows generalizations to other cones than $\mathbb R_+^d$. 
We keep the $\bm s,\nu$ notations in this article to simplify our formulas a little.

From now on and until the end of the paper, we consider $\xdata = \left(\xdataN\right) \in \mathbb R_+^{N\times d}$ an $N$-sample of independent and identically distributed (i.i.d.) random vectors, whose common distribution is the same as $\bm X$.

\begin{definition}[Monte-Carlo estimators]\label{ggchd:def:mc_estimators} For $\xdata \in \mathbb R_+^{N\times d}$ a dataset of $N$ independent and identically distributed (i.i.d.) samples of the $d$-variate random vector $\bm X$, define the estimators: 
\begin{align*}
  \hat{\bm\mu}(\xdata) &= \left(\widehat{\mu_{\bm k}}(\xdata)\right)_{\bm k \in I} = \left(\frac{1}{N}\sum_{i=1}^N \xdatai^{\bm k} e^{-\lvert \xdatai \rvert}\right)_{\bm k \in I}\\
  \hat{\bm\tau}(\xdata) &= \left(\widehat{\tau_{\bm k}}(\xdata)\right)_{\bm k \in I} = \bm B^{-1}\left(\hat{\bm\mu}(\xdata)\right).\\
  \hat{\bm a}(\xdata) &= \left(\widehat{a_{\bm k}}(\xdata)\right)_{\bm k \in I} =  \bm A \hat{\bm\mu}(\xdata) = \bm A \bm B(\hat{\bm\tau}(\xdata)),
\end{align*}
where $\hat{\bm \mu}(\xdata)$, $\hat{\bm \tau}(\xdata)$ and $\hat{\bm a}(\xdata)$ are all in $\mathbb R^{\lvert I \rvert}$. 
\end{definition}

\begin{remark}[Convergence and bias of $\hat{\bm\mu}(\xdata), \hat{\bm a}(\xdata)$ and $\hat{\bm\tau}(\xdata)$]\label{ggchd:rem:convergenc_eof_estimators}
All these estimators are convergent. Note that $\hat{\bm\mu}(\xdata)$ is the direct Monte-Carlo estimation of $\bm\mu$, and is therefore unbiased. Hence, $\hat{\bm a}(\xdata)$ is also an unbiased estimator of $\bm a$ (since the mapping between the two is linear), but $\hat{\bm\tau}(\xdata)$ is \emph{not} an unbiased estimator of $\bm \tau$.
Such unbiased estimators exist and are called $k$-statistics, see e.g.,~\cite{smith1995,smith2020,dinardo2009,nardo2011}, but this is not our focus here. 
\end{remark}

Following~\cite{laverny2021a}, we can therefore approximate the integrated square error between the densities through its projection in the truncated Laguerre basis, which provides the following empirical loss for a candidate Thorin measure $\nu$:  
$$\mathcal L(\xdata, \nu) = \lVert \hat{\bm a}(\xdata) - \bm A \bm B(\bm \tau(\nu))\rVert_2^2,$$
that can be minimized over an $n$-atomic Thorin measure $\nu$. In~\cite{laverny2021a}, it is shown that this loss benefits from a consistency property under the so-called well-behavior assumptions:

\begin{definition}[Well-behaved Gamma convolutions.] 
  The distribution $\mathcal G_{d}(\nu)$ is said to be well-behaved (w.b.) if and only if, denoting $S(\bm s)$ the ray in $\mathbb R_+^d$ that passes through $\bm 0$ and $\bm s$, 
  \begin{itemize}
    \item $\nu = \sum_{i=1}^n \alpha_i \delta_{\bm s_i}$ is finitely atomic with total mass $\lvert \nu \rvert = \sum_{i=1}^n \alpha_i \ge 1$.
    \item $\forall J \subseteq \{1,...,n\}$ such that $\nu\left(\BigCup\limits_{i \in J} S(\bm s_i)\right) > \nu\left(\BigCup\limits_{i \in I \setminus J} S(\bm s_i)\right)$, $\mathrm{Ker}\left(\left(\bm s_i\right)_{i\in J}\right) = \{\bm 0\}$.
  \end{itemize}
\end{definition}

By a slight abuse of language, we may talk about the well-behavior of the Gamma convolution itself or of its Thorin measure. Note that all univariate Gamma convolutions are well-behaved as soon as $\lvert\nu\rvert > 1$, and that well-behaved Gamma convolutions are absolutely continuous. Also note that we restricted to finitely atomic Thorin measures. The consistency property is then as follows:

\ifinthesis \restateconistencyggchd \else
\begin{theorem}[Consistency~\cite{laverny2021a}]\label{ggchd:thm:consistency_from_previous_paper}If $\xdata$ is drawn from a well-behaved density $f\in \mathcal G_{d}$,
any well-behaved estimator $\nu^\star$ such that $\mathcal L(\xdata, \nu^\star) \xrightarrow[N \to \infty]{a.s} 0$ ensures that
$$\lVert f - f_{\nu^\star} \rVert_2^2 \xrightarrow[\substack{N \to\infty\\I\to \mathbb N^d}]{a.s} 0.$$
\end{theorem}
\fi

Although it is consistent, the loss $\mathcal L$ is problematic because the function $\bm B$ is very combinatorial, and $\mathcal L(\xdata, \cdot)$ has a lot of local minimizers. Moreover, when the dimension $d$ increases, the number of coefficients $\lvert I\rvert$ that need to be matched must increase exponentially in the dimension to obtain correct results, and as soon as $d=4$ or $d=5$, the task of finding the global minimum of $\mathcal L(\xdata, .)$ becomes intractable. 

In this work, we propose to derive another loss, asymptotically equivalent to this one, and still consistent, that allows a better generalization to higher dimensional cases, avoiding the curse of dimensionality. We do this by using random projections. 

\section{Random projections and Grassmannian cubatures}\label{ggchd:sec:random_projs}
\ifinthesis\addtocontents{toc}{\protect\vskip-11pt}\fi

To simplify our notations and analysis, we assume an isotropic loss, using the increasing index set $I_m = \left\{\bm k \in \mathbb N^d, \lvert \bm k \rvert \le m\right\}$ for a given precision hyperparameter $m \in \mathbb N$ (where $\lvert \bm k \rvert = k_1 +...+k_d$). Its cardinal is $$D(m,d) = \sum\limits_{i=0}^m \binom{i+d-1}{d-1},$$ which increases exponentially in the dimension $d$ for a given precision $m$ of the approximation, as Table~\ref{ggchd:tab:D_m_d} clearly shows.

\tabletb
  \scriptsize
	\caption{\label{ggchd:tab:D_m_d}Number of coefficients $D(m,d)$.}
	\centering
	\begin{tabu} to \linewidth {>{\raggedright}l>{\raggedleft}X>{\raggedleft}X>{\raggedleft}X>{\raggedleft}X>{\raggedleft}X}
	\toprule
	\multicolumn{1}{c}{ } & \multicolumn{5}{c}{ } \\
  & $m=5$ & $m=10$ & $m=15$ & $m=20$ & $m=40$\\
	\cmidrule(l{3pt}r{3pt}){2-6} 
	
	\addlinespace[0.3em]
	$d=1$ & $6$   &$11$       &$16$        &$21$           &$41$\\
  $d=2$ & $21$&$66$      &$136$       &$231$          &$861$\\
  $d=3$ & $56$&$286$      &$816$      &$1771$        &$12341$\\
  $d=4$ & $126$&$1001$     &$3876$     &$10626$       &$135751$\\
  $d=5$ & $252$&$3003$    &$15504$     &$53130$      &$1221759$\\
  $d=6$ & $462$&$8008$    &$54264$    &$230230$      &$9366819$\\
  $d=7$ & $792$&$19448$   &$170544$    &$888030$     &$62891499$\\
  $d=8$ & $1287$&$43758$   &$490314$   &$3108105$    &$377348994$\\
  $d=9$ & $2002$&$92378$  &$1307504$  &$10015005$   &$2054455634$\\
  $d=10$ &$3003$&$184756$  &$3268760$  &$30045015$  &$10272278170$\\
  $d=20$ &$53130$&$30045015$  &$3247943160$  &$137846528820$  &$4191844505805495$\\
	\bottomrule
	\end{tabu}
\end{table}

As $D(m,d)$ is the number of Thorin moments that we need to compute at each iteration of a gradient descent on $\mathcal L$, we see that this loss is not usable for high dimensional problems. Computing more than $1000$ coefficients per gradient is already too much of a burden to be practical.
To solve this issue, we propose to use random projections. We now recall certain definitions and properties on the relationship that a measure has with its linear projections. 

\begin{definition}[Grassmannian] We denote $\mathcal E_k$ the $(k,d)$-Grassmannian, that is the set of $k$-dimensional linear subspaces $E$ of $\mathbb R^d$.
\end{definition}

\begin{definition}[Projection of a measure] Let $\nu$ be a measure on $\mathbb R^d$, and let $E \in \mathcal E_{k}$, $k < d$. We denote by $\nu_{\langle E\rangle}$ the projection of the measure $\nu$ onto this subspace, that is for $A \subseteq \mathbb R_+^k$, 
  $$\nu_{\langle E\rangle}(A) = \nu\left(\pi_{E}^{-1}(A)\right) = \nu\left(\left\{\bm x \in \mathbb R_+^d:\; \pi_{E}(\bm x) \in A\right\}\right),$$
  where $\pi_E$ is the linear projection onto the subspace $E$.
\end{definition}

As an exception, for any vector $\bm c \in \mathbb R^d$, we denote $\nu_{\langle \bm c\rangle}$ the projection onto the one-dimensional subspace generated by $\bm c$, identifying $\bm c$ and its span in the $d$-variate space, so that
$$\nu_{\langle \bm c\rangle}(A) = \nu\left(\left\{\bm x \in \mathbb R_+^d: \langle\bm c, \bm x\rangle \in A\right\}\right) \;\forall A \subseteq \mathbb R_+.$$

\begin{definition}[Agreement spaces] For two measures $\nu$ and $\mu$, define: 
  $$\mathcal E_k\left(\nu,\mu\right) = \left\{E \in \mathcal E_k, \nu_{\langle E\rangle} = \mu_{\langle E\rangle}\right\},$$
\end{definition}

Clearly, $\mathcal E_d\left(\nu,\mu\right) = \emptyset$ if $\nu \neq \mu$ and $\mathcal E_d\left(\nu,\mu\right) = \mathcal E_d$ if $\nu = \mu$. Note that $\mathcal E_d$ is a singleton, so $\mathcal E_d\left(\nu,\mu\right)$ has cardinal $1$ or $0$. The first result about projections comes from Cramer and Wold, back in 1936: 
\begin{theorem}[Cramer-Wold~\cite{cramer1936}] $\forall\, \nu_1,\nu_2 \in \mathcal M(\mathbb R^d),\; \nu_1 = \nu_2 \iff \mathcal E_1(\nu_1,\nu_2) = \mathcal E_1$.
\end{theorem}
Although this is interesting, more refined statements could be provided. In Appendix~\ref{ggchd:apx:projection}, we propose new results based on \cite{renyi1952,heppes1956} that deal with the particular case of projections of finitely atomic measures. 
These results are important as they ensure that the difference in value at each atom is bounded.
However, since we are fitting Thorin measures through a generalized moment problem (we try to match Laguerre coefficients), 
conditions on the first moments are sufficient for us. 
Weaker statements that only set the first moments can be obtained through the theory of Grassmannian cubatures, see~\cite{breger2017} and references therein. We describe a potential approach to these issues, inspired by~\cite{breger2017}.

We now denote by $\nu^{(\bm p)}$ the $\bm p^{\text{th}}$ raw moment of a measure $\nu$.

\begin{definition}[Moments agreement spaces] For two measures $\nu_1$ and $\nu_2$, $k\le d$ and an increasing set $I \subseteq \mathbb N^k$, define the moments' agreement space: 
  $$\mathcal E_{k,I}\left(\nu_1,\nu_2\right) = \left\{E \in \mathcal E_k,\; \nu_{1,\langle E\rangle}^{(\bm i)} = \nu_{2,\langle E\rangle}^{(\bm i)}\; \forall \bm i \in I\right\},$$
\end{definition}

Note that for all $I \subseteq \mathbb N^d$, $\mathcal E_{d,I}\left(\nu_1,\nu_2\right)$ has cardinal either $1$ or $0$, depending on the fact that $I$-moments of $\nu_1$ and $\nu_2$ match or not. 

Obviously, $\mathcal E_{k,I}\left(\nu_1,\nu_2\right) \subseteq \mathcal E_{k,J}\left(\nu_1,\nu_2\right)$ if $I \subseteq J$, and  $\mathcal E_{k,\mathbb N^d}\left(\nu_1,\nu_2\right) \subseteq \mathcal E_k\left(\nu_1,\nu_2\right)$. On the other hand $\mathcal E_{k,\mathbb N^d}\left(\nu_1,\nu_2\right) = \mathcal E_k\left(\nu_1,\nu_2\right)$ only if $\nu_2$ and $\nu_1$ are defined by their moments.

\begin{property}\label{ggchd:prop:linearisation_moments} Let $\mathcal C = \left\{\bm c_1,...,\bm c_{D(m,d)}\right\} \subseteq \mathcal E_1$ be such that the (square) matrix $$\mathcal P_{\mathcal C} = \left(\binom{m}{\bm k} \bm c_i^{\bm k}\right)_{\substack{i \le D(m,d)\\ \lvert\bm k\rvert \le m}}$$ is invertible. 

  \begin{enumerate}[label=(\roman*)]
    \item If for two measures $\nu_1$ and $\nu_2$, $\mathcal C \subseteq \mathcal E_{1,\{0,...,m\}}\left(\nu_1,\nu_2\right)$, then $\mathcal E_{d,I_m}\left(\nu_1,\nu_2\right) = \mathcal E_d$. That is, matching enough univariate moments is equivalent to matching multivariate moments.
    \item For any measure $\nu$, $$\left(\nu^{(\bm p)}\right)_{\bm p \in I_m} = \mathcal P_{\mathcal C}^{-1} \left(\nu_{\langle \bm c\rangle}^{(p)}\right)_{\substack{\bm c\in \mathcal C\\ p \le m}}.$$
    \item The matrix $\mathcal P_{\mathcal C}$ has bounded eigenvalues, and depends only on $\mathcal C$. 
  \end{enumerate}

\begin{proof}
We start by showing $(i)$. Since $\mathcal P_{\mathcal C}$ is invertible, its columns form a basis of a $D(m,d)$ dimensional space. 

  Denote by $\mathrm{Pol}_m\left(\mathbb R^d\right)$ the set of polynomials over $\mathbb R^d$ of degree at most $m$. A possible basis for this set is $$\left\{\bm x^{\bm \alpha} : \bm \alpha \in \mathbb N^d, \lvert \bm \alpha \rvert \le m \right\}.$$
  The dimension of $\mathrm{Pol}_m\left(\mathbb R^d\right)$ is therefore $D(m,d)$. Now, by the multinomial theorem, $$\langle \bm c_i, \bm x\rangle^m = \sum\limits_{\lvert \bm k\rvert \le m} \binom{m}{\bm k} \bm c_i^{\bm k} \bm x^{\bm k},$$
  Therefore, if $\mathcal P_{\mathcal C}$ is invertible, the polynomials $\left\{\langle \bm c, \bm x\rangle^m, \bm c \in \mathcal C\right\}$ form another basis of $\mathrm{Pol}_m\left(\mathbb R^d\right)$, and in particular the $d$-variate moments of $\nu_1$ and $\nu_2$ are a linear transformation of the univariate moments of $(\nu_{1,\langle c_i\rangle})_i$ and $(\nu_{2,\langle c_i\rangle})_i$, which implies $(i)$. 

  Moreover, $\mathcal P_{\mathcal C}$ is the transition matrix between the two bases, which implies $(ii)$. Finally, $\mathcal P_{\mathcal C}$ is of bounded operator norm as the linear projections of a given polynomial with finite coefficients cannot have infinite coefficients, which implies $(iii)$. 
\end{proof}
\end{property}

Of course, the hard part is finding such a set $\mathcal C$ making $\mathcal P_{\mathcal C}$ invertible. The following result is of particular interest in this matter:

\begin{theorem}[Random projections]\label{ggchd:thm:random_projs} If $\mathcal C$ is a set of $D(m,d)$ i.i.d. random vectors with $\mathcal U([0,1]^d)$ distributions, the random matrix $\mathcal P_{\mathcal C}$ is almost surely invertible.
  \begin{proof}
    Consider that $\bm c_1,...,\bm c_{D(m,d)}$ are i.i.d. random vectors uniform on $[0,1]^d$. 

    The matrix $\mathcal P_{\mathcal C}$ has a common factor $\left(\binom{m}{\bm k}\right)_{\lvert\bm k\rvert \le m}$ in each column. Hence, it is invertible if and only if $\mathcal Q_{\mathcal C}$ is, where: $$\mathcal Q_{\mathcal C} := \left(\bm c^{\bm k}\right)_{\substack{\bm c\in \mathcal C\\ \lvert\bm k\rvert \le m}}.$$

    Now, $\mathcal Q_{\mathcal C}$ is a generalized Vandermonde matrix. This kind of matrix is also called \emph{interpolation matrix} in a part of the literature, e.g., in~\cite{wendland2004}. 
    
    When $d=1$, $\mathcal Q_{\mathcal C}$ is a classical Vandermonde matrix, and therefore it is invertible if and only if the points $c_1,...,c_d \in [0,1]$ are distinct, which is almost sure.
    
    When $d > 1$, the matrix will be invertible if the set of points $\mathcal C$ is $\mathrm{Pol}_m(\mathbb R_+^d)$-unisolvant, as defined in \cite[Def. 2.6]{wendland2004}, that is if:
    
    $$\left\{P \in \mathrm{Pol}_m(\mathbb R_+^d):\; \forall \bm c \in \mathcal C,\;P(\bm c) = 0\right\} = \left\{0_{\mathrm{Pol}_m(\mathbb R_+^d)}\right\}.$$
    
    Therefore, for $\mathcal Q_{\mathcal C}$ not to be invertible, it would require that all points in $\mathcal C$ lie in the singular variety of a non-zero polynomial in $\mathrm{Pol}_m(\mathbb R_+^d)$. Such singular varieties are hypersurfaces of degree $m$. These singular varieties cannot have a positive Lebesgue measure in $[0,1]^{d}$, and therefore this almost surely never happens.
  \end{proof}
\end{theorem}

\begin{remark} Two small notes: 
  \begin{itemize}
    \item The condition that $\bm c_1,...,\bm c_{D(m,d)}$ are i.i.d. random vectors with $\mathcal U([0,1]^d)$ distributions is not even close to a necessary condition for the matrix $\mathcal P_{\mathcal C}$ to be invertible. Such necessary conditions are however hard to extract in all generality, but other nice particular examples\footnote{For example, a noteworthy $\mathrm{Pol}_m(\mathbb R_+^d)$-unisolvant set is $\mathcal C = \left\{\bm p \in \mathbb N^d : \lvert \bm p \rvert \le m\right\}$.} are given in \cite[Sec. 2]{wendland2004}. We do not know if there exists easily computatable better cubatures. 
    \item A (partial) generalization of Property~\ref{ggchd:prop:linearisation_moments} to dimensions of projections $k > 1$ might be possible, see, e.g., in~\cite{breger2017} recently. However, equivalent generalization of Theorem~\ref{ggchd:thm:random_projs} might not.
  \end{itemize}
\end{remark}

\section{Approximation of the Laguerre loss}\label{ggchd:sec:convex_approx}
\ifinthesis\addtocontents{toc}{\protect\vskip-11pt}\fi

In this section, we provide two successive approximations of the Laguerre loss that will ensure consistency, convergence and practical usability of the loss.
As we have already seen, the function $\bm B$ is bijective, but highly non-convex and has a complicated expression. The principal issues with the loss $\mathcal L$ are mostly due to the matrix $\bm A$ and the function $\bm B$: both are very combinatorial which might induce overflows if the implementation is not carefully done. They add a non-negligible computational cost to the minimization of $\mathcal L$, which has a myriad of local minimums. 

We propose to use a first order approximation of the function $\bm A\bm B$. Denote now $\nabla(\xdata)$ the Jacobian of the function $\bm \tau \mapsto \bm A \bm B (\bm \tau)$ taken at $\hat{\bm\tau}(\xdata)$. 

The function $\bm A \bm B$ can still be implemented efficiently, and its Jacobian $\nabla(\xdata)$ can be easily derived through symbolic computations. For lower dimensions, e.g., $d=1$, this can be done by hand to fasten the computations and such an algorithm is proposed in Appendix~\ref{ggchd:apx:nabla_AB}.

The function $\bm A\bm B$ being infinitely differentiable, it can be approximated by a first order Taylor expansion around $\hat{\bm\tau}(\xdata)$: 
\begin{align*}
  \bm A \bm B(\bm \tau (\nu)) &= \bm A \bm B(\hat{\bm \tau} (\xdata)) + \nabla(\xdata)\left(\bm \tau(\nu) - \hat{\bm\tau}(\xdata)\right) + \mathcal O\left(\lVert\bm \tau(\nu) - \hat{\bm\tau}(\xdata))\rVert_2^2\right)\\
  &= \hat{\bm a}(\xdata) + \nabla(\xdata)\left(\bm \tau(\nu) - \hat{\bm\tau}(\xdata)\right) + \mathcal O\left(\lVert\bm \tau(\nu) - \hat{\bm\tau}(\xdata))\rVert_2^2\right),
\end{align*}
which yields: 
\begin{align*}
  \mathcal L(\xdata,\nu) &=
\lVert \hat{\bm a}(\xdata) - \bm A\bm B\left(\bm \tau(\nu)\right)\rVert_2^2\\
      &= \lVert \nabla(\xdata) \left(\hat{\bm\tau}(\xdata) - \bm \tau(\nu)\right) \rVert_2^2 + \mathcal O\left(\lVert\bm \tau(\nu) - \hat{\bm\tau}(\xdata))\rVert_2^2\right)\\
      &= \lVert \hat{\bm\tau}(\xdata) - \bm \tau(\nu) \rVert_{\nabla(\xdata)}^2 + \mathcal O\left(\lVert\bm \tau(\nu) - \hat{\bm\tau}(\xdata))\rVert_2^2\right)
\end{align*}

We thus define, removing the higher order terms:
\[
  \widehat{\mathcal L}(\xdata, \nu)  = \lVert \hat{\bm\tau}(\xdata) - \bm \tau(\nu) \rVert_{\nabla(\xdata)}^2
\]

Note that the Jacobian $\nabla(\xdata)$ is a lower-triangular matrix and that $\nabla(\xdata)'\nabla(\xdata)$ is positive definite, which makes $\lVert \cdot \rVert_{\nabla(\xdata)}$ a proper norm.
We obtained a new loss $\widehat{\mathcal L}(\xdata, \cdot)$ to fit a $\mathcal G_{d}$ distribution on the data $\xdata$, from which we extract parameters by taking the argmin: $$\hat{\nu} = \argmin\limits_{\nu\, : \mathcal G_{d}(\nu)\, w.b.} \widehat{\mathcal L}(\xdata, \nu).$$

Computationally speaking, this is interesting since the complexity of the function $\bm A\bm B$ must only be dealt with once, at the beginning of the computations, and the result of its Jacobian $\nabla(\xdata)$ can be stored instead of recomputing the function $\bm A\bm B$ at each iteration. Moreover, this reduces the problem to a classic least-square generalized moment problem expression for the Thorin measure, with respect to the associated norm.  

Furthermore, this loss is still consistent as Property~\ref{ggchd:prop:consistency_widehat} shows. 
\ifinthesis \restateconssistencyofwidehatL \else
\begin{property}[Consistency of $\widehat{\mathcal L}$]\label{ggchd:prop:consistency_widehat} Let $\xdata$ be drawn from a well-behaved density $f\in \mathcal G_{d}$. The global minimizer
    $$\hat{\nu} = \argmin\limits_{\nu\, : \mathcal G_{d}(\nu)\, w.b.} \widehat{\mathcal L}(\xdata,\nu)$$
    ensures that $$\lVert f - f_{\hat{\nu}} \rVert_2^2 \xrightarrow[\substack{N \to\infty\\I\to \mathbb N^d}]{a.s} 0.$$
\end{property}
\fi 
\begin{proof}
      We already know by Theorem~\ref{ggchd:thm:consistency_from_previous_paper}, proved in~\cite{laverny2021a}, that any well-behaved estimator $\nu^*$ that ensures $$\mathcal L(\xdata,\nu^*) \xrightarrow[N \to \infty]{a.s.} 0$$ will have this consistency property.
      We therefore only need to show that $\mathcal L(\xdata, \hat{\nu}) \xrightarrow[N \to \infty]{a.s.} 0$.
      From the first-order Taylor expansion, we have: 
      \begin{align*}
        \mathcal L(\xdata, \hat{\nu}) &=
        \underbrace{\widehat{\mathcal L}(\xdata, \hat{\nu})}_{=: L_1} + \underbrace{\mathcal{O}\left(\lVert\hat{\bm\tau}(\xdata) - \bm\tau(\hat{\nu})\rVert_2^2\right)}_{=:L_2}.
      \end{align*}
      
      First, denote $\nu_{f}$ the Thorin measure of the true model $f$ and recall that it is well-behaved, and hence $\hat{\bm\tau}(x)$ and $\bm\tau(\hat{\nu})$ both converge to $\bm\tau(\nu_{f})$ almost surely: $\hat{\bm\tau}(\xdata)$ are Monte-Carlo estimators and therefore converge to $\bm\tau(\nu_{f})$ by Remark~\ref{ggchd:rem:convergenc_eof_estimators}, and $\bm\tau(\hat{\nu})$ minimizes a quadratic loss with target $\hat{\bm\tau}(\xdata)$. Hence, $L_2 \xrightarrow[N \to \infty]{} 0$.

      Then,  by the definition of $\hat{\nu}$ as the global minimizer of $\widehat{\mathcal L}$, $$\widehat{\mathcal L}(\xdata, \hat{\nu}) \le \widehat{\mathcal L}(\xdata, \nu_{f}) = \lVert \hat{\bm\tau}(\xdata) - \bm \tau(\nu_{f}) \rVert_{\nabla(\xdata)}^2,$$ which goes almost surely to zero since the Monte-Carlo estimator $\hat{\bm\tau}(\xdata) \xrightarrow[N \to \infty]{a.s.} \bm\tau(\nu_{f})$. 
      Hence, $\mathcal L(\xdata, \hat{\nu})$ converges to $0$, which concludes the argument. 
\end{proof}

Although the loss $\widetilde{\mathcal L}(\xdata, \cdot)$ is consistent,
  the number of coefficients $D(m,d) = \lvert I_{m} \rvert$ to compute at each iteration has not changed, and must still be exponentially increasing with the dimension.
Thus, we replace the loss $$\widehat{\mathcal L}(\xdata, \nu):= \lVert \hat{\bm\tau}(\xdata) - \bm \tau(\nu) \rVert_{\nabla(\xdata)}^2,$$
which is problematic computationally speaking, by another loss: 
\begin{align*}
  \widetilde{\mathcal L}(\xdata,\nu):&=  \int_{[0,1]^d} \widehat{\mathcal L}\left(\langle \bm c,\xdata\rangle, \nu_{\langle \bm c \rangle}\right) d\bm c\\
  &= \int_{[0,1]^d} \lVert \hat{\bm \tau}(\langle \bm c,\xdata\rangle) - \bm \tau(\nu_{\langle \bm c\rangle}) \rVert_{\nabla(\langle \bm c,\xdata\rangle)}^2 \,d\bm c,
\end{align*}
which we would minimize by stochastic gradient descent on $\bm c \sim \mathcal U([0,1]^d)$, drastically reducing the cost of a gradient iteration: this cost becomes linear in the precision $m$ and the dimension $d$ instead of exponential in $m,d$.
We now show the consistency of $\widetilde{\mathcal L}$.

\begin{theorem}[Consistency of $\widetilde{\mathcal L}$]\label{ggchd:thm:conssitency_of_widetile_L} Let $\xdata$ be drawn from a well-behaved density $f \sim \mathcal G_{d}(\nu)$. The global minimizer
  $$\tilde{\nu}:= \argmin\limits_{\nu\, : \mathcal G_{d}(\nu)\, w.b.} \widetilde{\mathcal L}(\xdata, \nu)$$
  ensures that $$\lVert f - f_{\tilde{\nu}} \rVert_2^2 \xrightarrow[\substack{N \to\infty\\I\to \mathbb N^d}]{a.s} 0.$$

  \begin{proof}

    According to Theorem~\ref{ggchd:thm:consistency_from_previous_paper}, we need to show that $\mathcal L(\xdata, \tilde{\nu}) \xrightarrow[N \to\infty]{a.s} 0$.

    We first show that $\widetilde{\mathcal L}(\xdata, \tilde{\nu}) \to 0$. For that, first note that $\tilde{\nu}$ is a global minimizer of $\widetilde{\mathcal L}(\xdata, \cdot)$ over the well-behaved parametrization. Denoting $\nu_{f}$ the (unknown) Thorin measure of the true model $f$ that generated the data, we suppose it to be well-behaved. It is therefore reachable by the minimization, and we have that: 

    $$\widetilde{\mathcal L}(\xdata, \tilde{\nu}) \le \widetilde{\mathcal L}(\xdata, \nu_{f}).$$

    Therefore, it suffices to show that $\widetilde{\mathcal L}(\xdata, \nu_{f}) \to 0$.

    For a given $\bm c \in [0,1]^d$, recall that $\mathcal L(\xdata, \nu_{f}) = \lVert \hat{\bm a}(\xdata) - \bm A\bm B\bm \tau(\nu_{f})\rVert_2^2$.
    Since the Laguerre functions verify $\lvert\phi_{\bm k}(\bm y)\rvert \le \sqrt{2}^d$ for any $\bm y \in \mathbb R_+^d$, we have 
    \begin{align*}
      \widehat{a_{\bm k}}(\xdata) &= \frac{1}{N}\sum_{i=1}^n \phi_{\bm k}(\xdatai) \le \sqrt{2}^d\\
      \bm A\bm B(\bm\tau(\nu_{f}))_{\bm k} &= \mathbb E\left(\phi_{\bm k}(\bm X)\right) \le \sqrt{2}^d
    \end{align*}
    which makes $\mathcal L(\xdata, \nu_{f})$ bounded, since there are a finite number $D(m,d)$ of coefficients in the sum. 
    
    Since the same argument can be done in a univariate case corresponding to the projection in a direction $\bm c$, $\mathcal L(\langle \bm c,\xdata\rangle, \nu_{f,\langle \bm c \rangle})$ is also bounded. Using the first order Taylor expansion, 
    $$\widehat{\mathcal L}(\langle \bm c,\xdata\rangle, \nu_{f,\langle \bm c \rangle}) = \mathcal L(\langle \bm c,\xdata\rangle, \nu_{f,\langle \bm c \rangle}) + \mathcal O\left(\lVert \hat{\bm a}(\langle \bm c, \xdata\rangle) - \bm A\bm B \bm \tau(\nu_{f,\langle \bm c \rangle})\rVert_2^2\right)$$ is also bounded. 
    Hence, almost all $\widehat{\mathcal L}(\langle \bm c,\xdata\rangle, \nu_{f,\langle \bm c \rangle})$ are bounded, positive, and go to zeros (by Monte-Carlo). Therefore, $\widetilde{\mathcal L}(\xdata, \nu_{f}) \to 0$ and hence $\widetilde{\mathcal L}(\xdata, \tilde\nu) \to 0$.

    We now show that $\mathcal L(\xdata, \tilde{\nu}) \xrightarrow[N \to\infty]{a.s} 0$.

    As we just showed, for almost all $\bm c \in [0,1]^d$, we have $\widehat{\mathcal L}\left(\langle \bm c,\xdata\rangle, \tilde{\nu}_{\langle \bm c \rangle}\right) \to 0$.
    Therefore, we can pick $D(m,d)$ elements in those $\bm c$ to form a set $\mathcal C$ making the matrix $\mathcal P_{\mathcal C}$ (almost surely) invertible by Theorem~\ref{ggchd:thm:random_projs}. Recalling that $\bm\tau$ are moments of a measure, we have:  
    \begin{align*}
      \widehat{\mathcal L}(\xdata, \tilde{\nu}) 
           &= \lVert \hat{\bm\tau}(\xdata) - \bm \tau(\tilde{\nu})\rVert_{\nabla(\xdata)}^2\\
           &= \lVert \mathcal P_{\mathcal C}^{-1}\left(\hat{\bm\tau}(\langle \bm c,\xdata\rangle)\right)_{\bm c \in \mathcal C} - \mathcal P_{\mathcal C}^{-1}\left(\bm \tau(\tilde{\nu}_{\langle \bm c \rangle})\right)_{\bm c \in \mathcal C}\rVert_{\nabla(\xdata)}^2\\
           &= \lVert \mathcal P_{\mathcal C}^{-1}\left(\hat{\bm\tau}(\langle \bm c,\xdata\rangle) - \bm \tau(\tilde{\nu}_{\langle \bm c \rangle}\rangle)\right)_{\bm c \in \mathcal C}\rVert_{\nabla(\xdata)}^2.\\
           &=\lVert\left(\hat{\bm\tau}(\langle \bm c,\xdata\rangle) - \bm \tau(\tilde{\nu}_{\langle \bm c \rangle}\rangle)\right)_{\bm c \in \mathcal C}\rVert_{\nabla(\xdata)\mathcal P_{\mathcal C}^{-1}}^2.
    \end{align*}

    We already showed that for almost all $\bm c$, $\widehat{\mathcal L}\left(\langle \bm c,\xdata\rangle, \tilde{\nu}_{\langle \bm c \rangle}\right) = \lVert \hat{\bm \tau}(\langle \bm c,\xdata\rangle) - \bm \tau(\tilde\nu_{\langle \bm c\rangle}) \rVert_{\nabla(\langle \bm c,\xdata\rangle)}^2 \to 0$, and hence $\hat{\bm\tau}(\langle \bm c, \xdata \rangle) - \bm\tau(\tilde\nu_{\langle \bm c \rangle}) \to 0$. Therefore, $\widehat{\mathcal L}(\xdata, \tilde{\nu}) \to 0$ since the linear map $\mathcal P_{\mathcal C}^{-1}$ is of bounded operator norm, and $\mathcal C$ is a finite set. Property~\ref{ggchd:prop:consistency_widehat} concludes.
  \end{proof}
\end{theorem}

The integral form of the loss $\widetilde{\mathcal L}$ allows stochastic gradient descent on $\bm c$: at each iteration, generate a random vector $\bm c$, and use the gradient  of $\widehat{\mathcal L}\left(\langle \bm c,\xdata\rangle, \nu_{\langle \bm c \rangle}\right)$. 
Although the stochasticity of the loss is on the direction of projection and not on the data, which is a little unconventional, we are not dealing with a true gradient descent in the sense that gradients of $\widehat{\mathcal L}\left(\langle \bm c,\xdata\rangle, \nu_{\langle \bm c \rangle}\right)$ are not subgradients of $\widetilde{\mathcal L}(\xdata, \nu)$.
Moreover, the global minimum does not intersect with the minimal varieties of the marginal losses, since the matrix $\nabla(\langle \bm c,\xdata\rangle)$ depends on the direction of projection.

Now that we have a tractable and consistent loss, and although the loss is not convex and might have local minimums, we show in the next section that the gradient flow associated to this loss converges to the global minimum under weak assumptions on the initialization. 

\section{Convergence of the gradient flow}\label{ggchd:sec:gradient_flows}
\ifinthesis\addtocontents{toc}{\protect\vskip-11pt}\fi

Although our loss $\widetilde{\mathcal L}$ has a precise structure, it is not convex. It has a myriad of local minimums. However, as we will now show, a direct and simple gradient descent on the atoms and weights of the Thorin measures with a large enough number of proposed atoms, spread densely enough in the space, will converge to the global minimum under some very minor assumptions on the regularity of the true Thorin measure.

To postpone the question of which algorithm is used to minimize our objective function $\widetilde{\mathcal L}$, we consider the framework of gradient flows on the space $\mathcal M_+(\mathbb R_+^d)$ of positive Thorin measures. A numerical solver that uses a gradient descent on the loss, or a stochastic gradient descent, will in fact approximate the gradient flow of the loss (the sequence of iterates converges to the gradient flow when the step size of the gradient descent goes to zero). 
After some definitions, we show that our settings are covered by the analysis from Chizat~\cite{chizat2021} (see also~\cite{chizat2018}), which ensures global convergence of the gradient flow. 

From now on, let $\mathcal F = L_2([0,1]^d,\mathbb R_+^m)$ be the space of square integrable functions from $[0,1]^d$ to $\mathbb R_+^m$.

We consider the functional $\psi$ defined as:
\begin{equation*}
  \psi\colon \left\{\begin{aligned}
    \mathcal \mathbb R_+^d   &\longrightarrow \mathcal F\\
    \bm s       &\longmapsto \psi(\bm s)\colon \left\{\begin{aligned}
      [0,1]^d   &\longrightarrow \mathbb R_+^m\\
        \bm c      &\longmapsto\bm \tau\left(\delta_{\bm s, \langle \bm c \rangle}\right).
     \end{aligned}\right.
   \end{aligned}\right.
\end{equation*}
Note that $\delta_{\bm s, \langle \bm c \rangle} = \delta_{\langle\bm s, \bm c \rangle}$ is the univariate Dirac measure in $\langle \bm s,\bm c\rangle$.
With a slight abuse of notations, we also denote 
\begin{equation*}
  \psi(\nu) = \int \psi(\bm s) \nu(d\bm s) \colon \left\{\begin{aligned}
    [0,1]^d   &\longrightarrow \mathbb R_+^m\\
      \bm c      &\longmapsto \bm\tau(\nu_{\langle \bm c\rangle}) = \int_{\mathbb R_+^d} \bm \tau\left(\delta_{\bm s, \langle \bm c \rangle}\right)\nu(d\bm s)
   \end{aligned}\right. \;\in \mathcal F.
\end{equation*}

We also consider, as our target, 
\begin{equation*}
  \psi_{\xdata} \colon\left\{\begin{aligned}
    [0,1]^d   &\longrightarrow \mathbb R_+^m\\
      \bm c      &\longmapsto \hat{\bm \tau}(\langle \bm c, \xdata \rangle)
   \end{aligned}\right. \;\in \mathcal F.
\end{equation*}

We endow $\mathcal F$ with the norm $\lVert \cdot \rVert_{\xdata}$, which depends on the random sample $\xdata$,
defined (since all $\nabla(\langle \bm c, \xdata\rangle)$ are positive definite) by: 
$$\forall \psi \in \mathcal F,\;\lVert \psi \rVert_{\xdata}^2 = \int_{[0,1]^d} \rVert \psi(\bm c)\lVert_{\nabla(\langle \bm c, \xdata\rangle)}^2 d\bm c,$$
where the norms $\lVert \cdot \rVert_{\nabla(\cdot)}^2$ were previously defined in Section~\ref{ggchd:sec:convex_approx}. The scalar product associated with the norm $\lVert \cdot \rVert_{\xdata}$ is denoted $\langle \cdot, \cdot \rangle_{\xdata}$. The corresponding bilinear operator can be seen as a blockwise diagonal operator with diagonal blocks $\nabla(\langle \bm c, \xdata \rangle)'\nabla(\langle \bm c, \xdata \rangle)$ for all $\bm c \in [0,1]^d$. This norm makes $\mathcal F$ a Hilbert space, and it furthermore allows to write our loss as: 
$$\widetilde{\mathcal{L}}(\xdata, \nu) = \lVert \psi_{\xdata} - \psi(\nu)\rVert_{\xdata}^2,$$
which shows that our loss in convex in $\psi(\nu)$, but clearly not convex in the atoms and weights of an $n$-atomic measure $\nu$.

The idea behind the conic particle gradient descent we propose is to parametrize the candidate measure $\nu$ as $$\nu = \sum_{i=1}^n p_i^2 \delta_{\bm s_i}$$ and to perform direct gradient descent on $\left(\bm p, \bm s\right) \mapsto \widetilde{\mathcal L}(\xdata,\sum_{i=1}^n p_i^2 \delta_{\bm s_i})$. If the chosen number of atoms $n$ is large enough, and if the initial $\bm p, \bm s$ are spread enough in the space, then the gradient descent will converge to a Thorin measure that achieves the global minimum. 

If we endow the space of measures $\mathcal M(\mathbb R_+^d)$ with the total variation norm, the objective $\widetilde{\mathcal{L}}(\xdata, \cdot)$ is Fréchet differentiable, and we can compute its differential. 

\begin{property} The Fréchet differential of $\widetilde{\mathcal{L}}(\xdata, \cdot)$ at $\nu \in \mathcal M(\mathbb R_+^d)$ is represented by the function 
  \begin{equation*}
    \widetilde{\mathcal L}_{\nu}' \colon\left\{\begin{aligned}
      \mathbb R_+^d  &\longrightarrow \mathbb R\\
        \bm s      &\longmapsto 2\langle \psi_{\xdata} - \psi(\nu), \psi(\bm s) \rangle_{\xdata}.
     \end{aligned}\right.
  \end{equation*}
  
  \begin{proof}
    The differential of $\widetilde{\mathcal{L}}(\xdata, \cdot)$ at $\nu \in \mathcal M(\mathbb R_+^d)$ is defined as the function $\widetilde{\mathcal L}_{\nu}': \mathbb R_+^d \to \mathbb R$ such that for any measure $\sigma \in \mathcal M(\mathbb R_+^d)$, 
    $$\left. \frac{\partial}{\partial \epsilon} \widetilde{\mathcal{L}}(\xdata,\nu + \epsilon \sigma) \right\lvert_{\epsilon = 0} = \int \widetilde{\mathcal L}_{\nu}'(\bm s) d\sigma(\bm s)$$
    
    Which yields:
  \begin{align*}
    \left. \frac{\partial}{\partial \epsilon} \widetilde{\mathcal{L}}(\xdata,\nu + \epsilon \sigma) \right\lvert_{\epsilon = 0}
    &= \left. \frac{\partial}{\partial \epsilon} \lVert \psi_{\xdata} - \psi(\nu + \epsilon\sigma)\rVert_{\xdata}^2 \right\lvert_{\epsilon = 0} \\
    &= \left. 2\langle \psi_{\xdata} - \psi(\nu) - \epsilon\psi(\sigma),\psi(\sigma)\rangle_{\xdata} \right\lvert_{\epsilon = 0} \\
    &= 2\langle \psi_{\xdata} - \psi(\nu),\psi(\sigma)\rangle_{\xdata} \\
    &= 2\langle \psi_{\xdata} - \psi(\nu),\int \psi(\bm s) \;\sigma(d\bm s)\rangle_{\xdata} \\
    &= \int 2\langle \psi_{\xdata} - \psi(\nu), \psi(\bm s) \rangle_{\xdata}\; \sigma(d\bm s).
  \end{align*}
  \end{proof}
\end{property}

Following~\cite{chizat2021}, we lift the problem in the \emph{Wasserstein space} of probability measures, by considering probability measures with atoms in $d+1$ dimensions, the first one corresponding to our weights. For any atom $(p,\bm s)$, we define the Riemannian inner product in this space by:
\begin{equation}\label{ggchd:eq:metric}
  \langle (p_1,\bm s_1), (p_2, \bm s_2)\rangle_{(p,\bm s)} = p_1p_2 + p^2\langle \bm s_1, \bm s_2\rangle.
\end{equation}

Recall that the Wasserstein distance between two probability measures is defined as $$W_{p}(\mu_1,\mu_2) = \left(\min\limits_{\gamma \in \Pi(\mu_1,\mu_2)} \int \mathrm{dist}(\bm x,\bm y)^p \gamma(\partial \bm x,\partial \bm y) \right)^{\frac{1}{p}},$$
where the distance $\mathrm{dist}(\bm x, \bm y)$ is associated to the cone metric from Equation~\eqref{ggchd:eq:metric}. Moreover, we define the $\max$-Wasserstein distance $W_{\infty}$ as the limit of $W_{p}$ when $p$ goes to infinity, or more directly as $$W_{\infty}(\mu_1,\mu_2) = \inf\limits_{\gamma \in \Pi(\mu_1,\mu_2)} \max_{(\bm x, \bm y) \in \mathrm{Sup}(\gamma)} \mathrm{dist}(\bm x,\bm y).$$

\begin{property}\label{ggchd:prop:gradient}
  The gradient $\nabla \widetilde{\mathcal L}$ of $\widetilde{\mathcal{L}}(\xdata, \cdot)$ w.r.t. the conic metric from Equation~\eqref{ggchd:eq:metric}, for $\nu = \sum\limits_{i=1}^n p_i^2 \delta_{\bm s_i}$, is defined componentwise as:
  \begin{align}
    \frac{\partial \widetilde{\mathcal{L}}(\xdata, \nu)}{\partial p_i} &= 2p_i \widetilde{\mathcal L}_\nu'(\bm s_i) = 4p_i \langle \psi_{\xdata} - \psi(\nu), \psi(\bm s_i) \rangle_{\xdata}\\
    \frac{\nabla \widetilde{\mathcal{L}}(\xdata, \nu)}{\partial \bm s_{i}} &= \nabla\widetilde{\mathcal L}_\nu'(\bm s_i) = 2\langle \psi_{\xdata} - \psi(\nu),  \frac{\nabla\psi(\bm s_i)}{\partial \bm s_{i}} \rangle_{\xdata},
  \end{align}
  where the gradient $\frac{\nabla\psi(\bm s)}{\partial \bm s}$ of $\psi$ has components $\frac{\partial \psi(\bm s)}{\partial s_{j}} = \left(\frac{\partial \bm \tau\left(\delta_{\langle \bm s, \bm c \rangle}\right)}{\partial s_{j}}\right)_{\bm c \in [0,1]^d},$
whose elements can be easily derived from the expression of $\bm\tau$ given in Property~\ref{ggchd:prop:thorin_moments}:
\begin{align*}
  \frac{\partial \bm \tau_{0}(\delta_{\langle \bm c,\bm s\rangle})}{\partial s_j} 
    = -\frac{\partial \ln\left(1 + \langle \bm c,\bm s\rangle\right)}{\partial s_j}
    &= -\frac{s_j}{1 + \langle \bm c,\bm s\rangle}, \text{ and}\\
  \frac{\partial \bm \tau_{k}(\delta_{\langle \bm c,\bm s\rangle})}{\partial s_j} 
    = \frac{\partial \langle \bm c,\bm s\rangle^{k}\left(1 + \langle \bm c,\bm s\rangle)\right)^{-k}}{\partial s_j}
    &= \frac{k c_j \langle \bm c, \bm s \rangle^{k-1}}{\left(1 +\langle \bm c, \bm s \rangle\right)^{k+1}},
  \text{ for } k = 1,...,m.\\ 
\end{align*}
\end{property}

\begin{remark}[Metric choice] The choice of metric on the lifted Wasserstein space decides how to combine the \emph{vertical} and \emph{spatial} components of the gradient. The choice we made, which allows further results from Theorem~\ref{ggchd:thm:convergence_of_gradient_flow} to apply, is to decouple the two and have the spatial components of the gradient be independent of the vertical ones. 
\end{remark}

Now that the gradient $\nabla \widetilde{\mathcal L}$ is well-defined, we consider its flow as given by Definition~\ref{ggchd:def:gradient_flow}:

\ifinthesis \restategradientflow \else
\begin{definition}[Gradient flow\cite{ambrosio2005}]\label{ggchd:def:gradient_flow} A gradient flow of $\widetilde{\mathcal{L}}(\xdata, \cdot)$ is an absolutely continuous curve $\left(\nu_t\right)_{t \ge 0}$ in the space $\mathcal M(\mathbb R_+^d)$ that satisfies 
  $$\frac{\partial}{\partial t}\nu_t = - \nabla \widetilde{\mathcal L}(\xdata, \nu_t)$$ 
\end{definition}
\fi

The gradient flow $\left(\nu_t\right)_{t\ge 0}$, with initialization $\nu_0$, converges to a certain Thorin measure $\nu_\infty$, local minimizer of $\widetilde{\mathcal L}$. Theorem~\ref{ggchd:thm:convergence_of_gradient_flow}, adapted from~\cite{chizat2021}, concludes this analysis by ensuring that $\nu_\infty$ is a global minimizer, under certain regularity conditions on the initialization $\nu_0$. 

\ifinthesis \restateconvergenceflow \else
\begin{theorem}[Global convergence of the gradient flow]\label{ggchd:thm:convergence_of_gradient_flow}
For $\rho \in \mathcal M_+(\mathbb R_+^d)$ an absolutely continuous reference measure such that $\log\rho$'s density is Lipschitz, for an initial measure $\nu_0 \in \mathcal M_+(\mathbb R_+^d)$, there exists a constant $C$, dependent on the characteristics of the problem, such that if $W_{\infty}(\nu_0,\rho) \le C,$ $$\exists\, \nu_{\infty} \in \argmin\limits_{\nu \in \mathcal M_+(\mathbb R_+^d)} \widetilde{\mathcal L}(\xdata, \nu) \text{ such that } W_{\infty}(\nu_t,\nu_{\infty}) \xrightarrow[t \to \infty]{} 0.$$

Furthermore, when $\nu_0 = \rho$, we achieve a precision $\epsilon$, i.e., $W_{\infty}(\nu_t,\nu_{\infty}) \le \epsilon$, provided the number of iteration is $t = \mathcal O\left(-\log(\epsilon)\right)$.
\end{theorem}
\fi
\begin{proof}
    This is a direct application of Chizat's result~\cite[Th. 4.1]{chizat2021}. See also~\cite{chizat2018} and~\cite[Apx. A.3]{decastro2021} for the complexity.
      %
      %
\end{proof}

For general accelerated convex methods, $W_\infty(\nu_t,\nu_\infty)\le \epsilon$ is achieved for $t = \mathcal O\left(\epsilon^{-1/d}\right)$ only~\cite{decastro2021}.

Although Chizat's theorem applies in more general settings than ours, and gives more precise results, the constants are not of so much importance to us. What matters is that, when initialized with enough atoms, picked densely enough in the space (so that the atomic measure $\nu_0$ is \emph{close enough}, in the Wasserstein $W_\infty$ sense, to a $\log$-Lipschitz measure $\rho$), the (maybe stochastic) gradient descent on $\widetilde{\mathcal{L}}(\xdata, \cdot)$ converges to the global minimizer of the problem.
A good initialization is therefore, e.g., to take all components of $\bm p, \bm s$ independent centered Gaussian (for $\bm p$) and independent chi squared (for $\bm s$, which must stay positive). A very bad initialization will be to take all atoms in the same spot, e.g., $\bm p \equiv 1$, $\bm s \equiv 1$. Finally, initializing atoms clustered in a small space is also not a good idea. 

Chizat's analysis gives a nice feature to the final structure of the algorithm: the possibility to include a Lasso constraint.

\begin{remark}[Lasso constraint]\label{ggchd:rem:lasso} If one is concerned with the sparsity of the produced Thorin measure, an additional total variation regularization can be added by considering a hyperparameter $\lambda > 0$ and minimizing the loss 
  $$\widetilde{\mathcal{L}}(\xdata, \nu) + \lambda \lvert \nu \rvert.$$
  Results from Theorem~\ref{ggchd:thm:convergence_of_gradient_flow} still hold, and any wanted hyperparameter searching method (with its associated cost) could be used to find $\lambda$.
\end{remark}

\section{Performance investigation} \label{ggchd:sec:investigation}
\ifinthesis\addtocontents{toc}{\protect\vskip-11pt}\fi

We will minimize $\widetilde{\mathcal L}$ by stochastic gradient descent on $\bm c$.
Indeed, as the dimension increases, the possibility to explore at least $D(m,d)$ different directions of projections vanishes, and our numerical capacity to compute a gradient of $\mathcal L$ or $\widehat{\mathcal L}$ vanishes with it. We initialize the Thorin measure by setting its atoms and weights as random Gaussian noise, to spread them enough, according to the analysis from Section~\ref{ggchd:sec:gradient_flows}. The computation of gradients is discussed thoroughly in Appendix~\ref{ggchd:apx:nabla_AB}. To improve stability of the Laguerre expansions, we also preprocess the datasets so that each marginal has unit variance. Our algorithm is then the following: 

\algotb
	\SetAlgoLined
	\KwIn{A dataset $\xdata \in \mathbb{R}^{N\times d}$, a number of Gammas $n \in \mathbb N$, a precision parameter $m \in \mathbb N$, a number of iterations $T \in \mathbb N$, and a learning rate $\eta \in \mathbb R_+$}
	\KwResult{A Thorin measure $\nu_T$ that approximates the dataset $\xdata$ as a multivariate Gamma convolution.}
  Estimate standard deviations $\sigma_i = \mathrm{std}(\xdatai)$ for all $i \in 1,...d$, and standardize the marginals by dividing $\xdatai$ by $\sigma_i$.\\
  Initialize a measure $\nu_0 \in \mathcal M_{+}(\mathbb R_+^d)$ with $n$ atoms and corresponding weights, chosen randomly through Gaussian noise.\\
  \ForEach{$t \in 0,...,T-1$}{
    Choose a random direction $\bm c \in [0,1]^d$. \\
    Compute the gradient $\bm g$ of $\widehat{\mathcal L}(\langle \bm c, \xdata\rangle, \nu_{t,\langle \bm c\rangle})$ with respect to $\nu_t$ according to Property~\ref{ggchd:prop:gradient} and Appendix~\ref{ggchd:apx:nabla_AB}. \\
    Let $\nu_{t+1} = \nu_t - \eta\bm g$
	}
  Rescale $\nu_T$ by $\sigma_1,...\sigma_d$.\\
	Return $\nu_T$
	\caption{Estimation of Thorin measures via stochastic gradient descent on $\widetilde{\mathcal L}$.}
  \label{ggchd:algo:fitting_algo}
\end{algorithm}

Note that the expectation of the stochastic gradient is indeed the gradient of the global loss (by the very definition of the global loss $\widetilde{\mathcal L}$). Moreover, these stochastic gradients are all independents knowing $\xdata$.

In Algorithm~\ref{ggchd:algo:fitting_algo}, the number of iterations $T$ is supposed to be chosen on a case-by-case basis: by default, our implementation runs until (about) $\lvert \bm g \rvert \le 1e-16.$\footnote{for all $X \in\mathbb N$, $1e-X$ is a notation that represents the number $10^{-X}$. $1e-16$ is approximately the \texttt{Float64} machine epsilon, according to the IEEE754-2008 standard.}

Moreover, $\eta$ is simply a fixed step size (or \emph{learning rate}) for the gradient descent. However, in practice, as our parameters are very different things (atoms and weights of the Thorin measure are of different nature), we chose to use different learning rates per parameters, and also through time via momentum.
The standard algorithm Adam~\cite{kingma2014,dozat2016} seems to be a good fit for our problems. Indeed, the adaptive nature of Adam learning rates (one per parameter) is a good numerical solution to the problem of different vertical and spatial scaling of the gradients (on one side) and to different learning rate along the fitting process (on the other side).

To compare multivariate distributions to assess our results, we use non-distribution-free versions of the Kolmogorov-Smirnov and Cramer-Von~Mises tests which are detailed in depth in Algorithms~\ref{ggchd:algo:benchmark} and~\ref{ggchd:algo:ks_cm_stats}, in Appendix~\ref{ggchd:apx:KSCvM}. By construction, these tests are calibrated correctly: their p values are uniform if and only if the distributions are the same. Therefore, we plot histograms of (resamples of) these p values and check for approximate uniformity as in Figure~\ref{ggchd:fig:ks_test_histogram}. However, these tests are limited to low dimensional cases, as the task of computing an empirical version of the distribution of the statistics under the null hypothesis has a complexity that is exponential in the dimension.

We start by a simple example in four dimensions.

\subsection{A distribution with functional structure.}

The following example deals with a particular dependence structure that we generated with a certain functional, as Dataset~\ref{ggchd:data:WeirdExample} describes. 

\ifinthesis \restateweirddataset \else
\begin{dataset}[Functional dependence structure]\label{ggchd:data:WeirdExample} For $U_1,U_2,U_3$ independent $\mathcal U([0,1])$ distributions and  $Y_1,...,Y_4$ independent log-Normal(0,1) distributions, we let the random vector $\bm X$ be defined by: 
  $$\bm X = \left(Y_1,\, Y_2 + U_1 Y_1^2,\, Y_3 + U_3 Y_1,\, Y_4 + Y_1^{1 + \frac{U_3}{3}}\right),$$
  We simulate a dataset $\bm x \in \mathbb R^{10000 \times 4}$ of $N = 10000$ i.i.d. samples from $\bm X$.
\end{dataset}
\fi
This produces a weird dataset that is typically hard to fit with classical parametric models such as (hierarchical) Archimedean copulas. However, as Figure~\ref{ggchd:fig:WeirdExample_summary} will show, this distribution is not too far from the $\mathcal G_d$ class.

We launch the optimization of $\widetilde{\mathcal L}$ on Dataset~\ref{ggchd:data:WeirdExample}, with a starting parametrization with $n=100$ Thorin atoms. We let the optimizer run until convergence, which always appends since our learning rates are decreasing. Note that we have no way of checking the global optimality of our solutions, but we can extract some insights from Chizat's results: if several atoms end up in at the same spot, yielding a final measure with much less than $n=100$ atoms, this is a good sign for convergence. 
In our simulations this behavior is common: the optimizer tends to set Thorin weights to zero if the atom is not needed, or to merge atoms that are close by. 

\figuretb
  \centering
  \includegraphics[width=\linewidth]{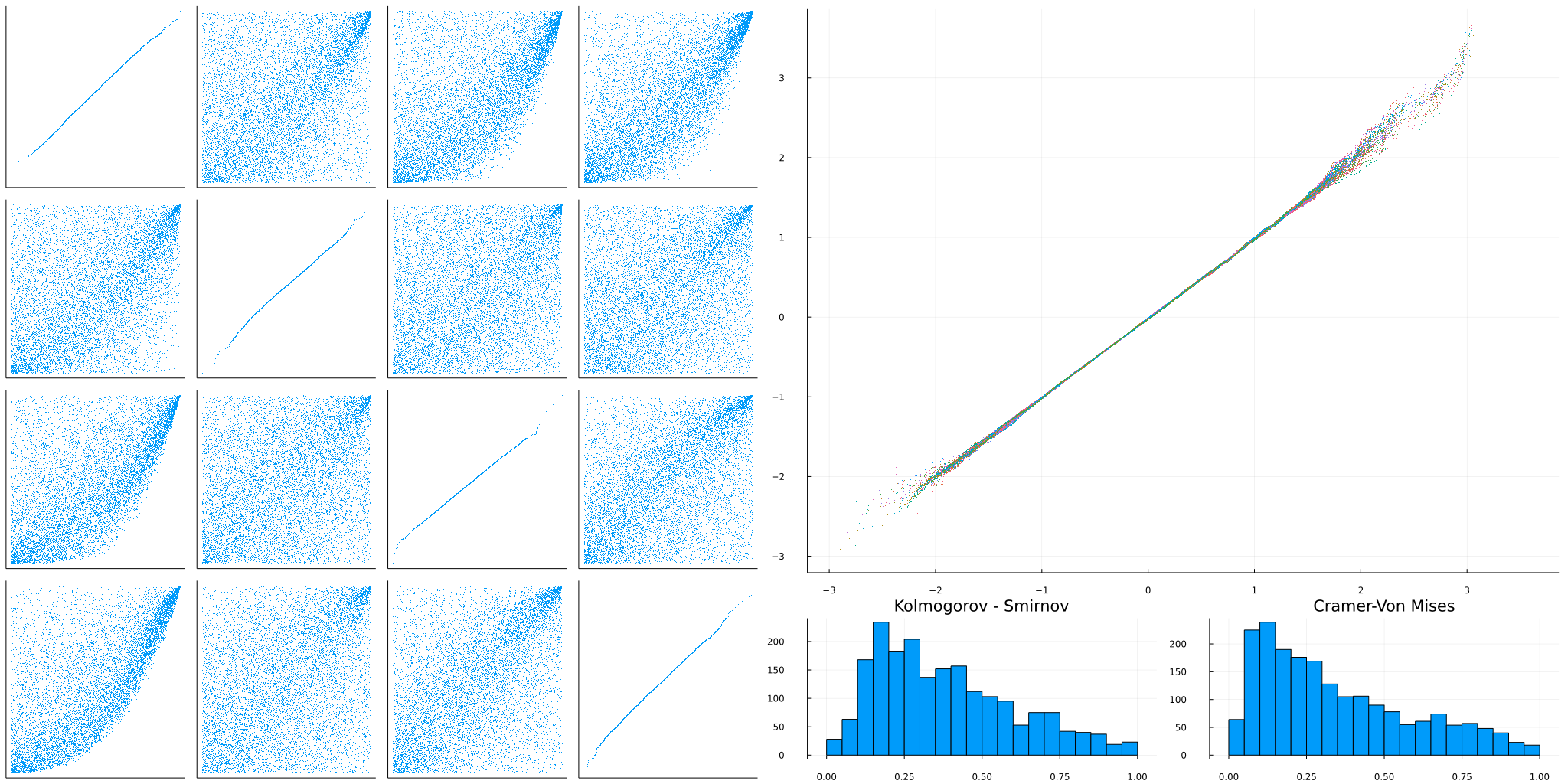}
  \caption{Results on Dataset~\ref{ggchd:data:WeirdExample}. On the left, a pairs plot representing the four dimensions: the bottom left triangle represents the pseudo-data from the original dataset, while the top triangle represent the pseudo-data from our estimator, and the middle line represents quantile-quantile plots of the marginals. Top right, quantile-quantile plots of $\langle \bm c_i, \hat{\bm X}\rangle$ on randomly chosen directions $\bm c_1,...,\bm c_{50}$. Bottom right, histograms of approximated Kolmogorov-Smirnov and Cramer-Von~Mises p values of our estimation.}
  \label{ggchd:fig:WeirdExample_summary}
\end{figure}

The estimator depicted in Figure~\ref{ggchd:fig:WeirdExample_summary} has a Thorin measure with $n=37$ atoms, although we started the estimation with $n=100$. Most of these atoms have a very low contribution to the mean of the marginals however (very small weight $\alpha$ and/or very small scale $s$), as Table~\ref{ggchd:tab:thorin_measure_WeirdExample} depicts.

\tabletb
  \scriptsize
	\caption{\label{ggchd:tab:thorin_measure_WeirdExample}Thorin measure obtained on Dataset~\ref{ggchd:data:WeirdExample}}
	\centering
	\begin{tabu} to \linewidth {>{\raggedleft}l>{\raggedleft}X>{\raggedleft}X>{\raggedleft}X>{\raggedleft}X}
	\toprule
  $\alpha$    & $s_1$       & $s_2$       & $s_3$       & $s_4$\\
	\midrule
	\addlinespace[0.3em]
  $4.855e-16$ & $8.627e-8$  &     --      &     --      & $0.914$\\
  $3.277e-15$ & $291.3$     & $2.095e-10$ & $1.043e-8$  & $2.513e-7$\\
  $2.921e-12$ &     --      &     --      & $7.433e-5$  & $0.9138$\\
  $3.756e-12$ & $279.4$     & $9.078e-5$  & $0.9684$    & $2.558e-11$\\
  $4.250e-12$ &     --      &     --      & $1.171e-10$ & $40.67$\\
  \addlinespace[0.3em]
  $1.528e-11$ & $0.09957$   &     --      &     --      &      --      \\
  $5.249e-11$ & $7.407e-13$ & $5.33e-15$  & $4.872e-11$ & $0.9136$\\
  $6.004e-11$ &     --      &     --      &     --      & $0.9138$\\
  $6.761e-11$ & $5.164e-16$ & $26.9$      &     --      &      --      \\
  $1.337e-9$  & $4.533e-8$  &     --      &     --      & $0.9139$\\
  \addlinespace[0.3em]
  $2.193e-8$  &     --      & $2.394e-13$ &     --      & $0.9137$\\
  $3.860e-8$  &     --      & $2.261e-15$ &     --      & $0.9134$\\
  $4.031e-8$  &     --      & $0.1374$    & $3.583e-14$ & $0.1533$\\
  $1.036e-7$  & $6.337e-13$ & $1.053e-8$  &     --      & $0.9138$\\
  $1.368e-7$  &     --      &     --      &     --      & $0.9132$\\
  \addlinespace[0.3em]
  $2.826e-7$  & $1.218e-14$ & $1.669e-13$ &     --      & $0.9138$\\
  $5.642e-7$  &     --      &     --      &     --      & $40.69$\\
  $1.089e-6$  &     --      &     --      &     --      & $0.9144$\\
  $1.473e-6$  &     --      & $2.607e-10$ & $0.000824$  & $40.66$\\
  $1.905e-6$  &     --      &     --      &     --      & $0.9138$\\
  \addlinespace[0.3em]
  $1.114e-5$  & $5.208e-15$ &     --      & $0.0003181$ & $0.914$\\
  $7.958e-5$  & $5.706e-15$ &     --      &     --      & $0.9137$\\
  $0.0001191$ &     --      &     --      &     --      & $0.9137$\\
  $0.0001421$ &     --      &     --      & $9.108e-12$ & $0.9137$\\
  $0.005818$  &     --      & $0.322$     & $7.331$     & $0.09808$\\
  \addlinespace[0.3em]
  $0.01199$   & $4.534$     & $7.889$     & $2.563$     & $7.845$\\
  $0.1074$    & $1.238$     &     --      & $0.9892$    & $1.092$\\
  $0.1865$    &     --      & $1.432e-15$ &     --      & $0.9138$\\
  $0.2146$    & $1.316$     & $0.8674$    & $0.8449$    & $0.8522$\\
  $0.2167$    & $3.835e-8$  & $1.277$     & $0.03307$   & $1.447e-8$\\
  \addlinespace[0.3em]
  $0.2769$    & $7.69e-14$  &     --      & $1.006$     &      --      \\
  $0.3426$    & $0.09957$   &     --      & $1.073e-7$  & $0.0004082$\\
  $0.4323$    & $0.1612$    & $0.1614$    & $0.08263$   & $0.04896$\\
  $0.5334$    & $0.2333$    & $5.798e-8$  & $0.2282$    & $0.07975$\\
  $1.539$     & $2.931e-9$  & $0.1434$    &     --      &      --      \\
  \addlinespace[0.3em]
  $1.745$     & $0.01782$   &     --      & $0.01409$   & $0.1335$\\
  $3.797$     &     --      &     --      & $0.09868$   &      --      \\
	\bottomrule
	\end{tabu}
\end{table}

Several observations can be made on the results from Table~\ref{ggchd:tab:thorin_measure_WeirdExample}. The first one is that most of the Thorin weights are concentrated in only a few atoms, and many atoms have weights and/or values that are very low, the minimum being $\sim 1e-16$. This is because we did not set a threshold beside the \texttt{Float64} limits imposed by our computations. One could argue that Thorin weights below a certain value are just noise and should be removed. If we set a zero-threshold at $1e-05$ on both the weights and components of atoms, the remaining atoms are only the $17$ last lines, with even more zeros in the scales.

This 'thresholded' estimator gives results depicted on Figure~\ref{ggchd:fig:WeirdExample_threasold_summary}. We see that this threshold fixed a little the problem in the projections tails, and we obtained better results on the Cramer-Von~Mises and the Kolmogorov-Smirnov tests, without changing the pairs plot of the pseudo-observations. For this reason, we included a thresholding possibility in the code, with a user defined minimal Thorin weight. We do not provide analysis on these thresholds, but a link to Remark~\ref{ggchd:rem:lasso} could maybe be drawn.

\figuretb
  \centering
  \includegraphics[width=\linewidth]{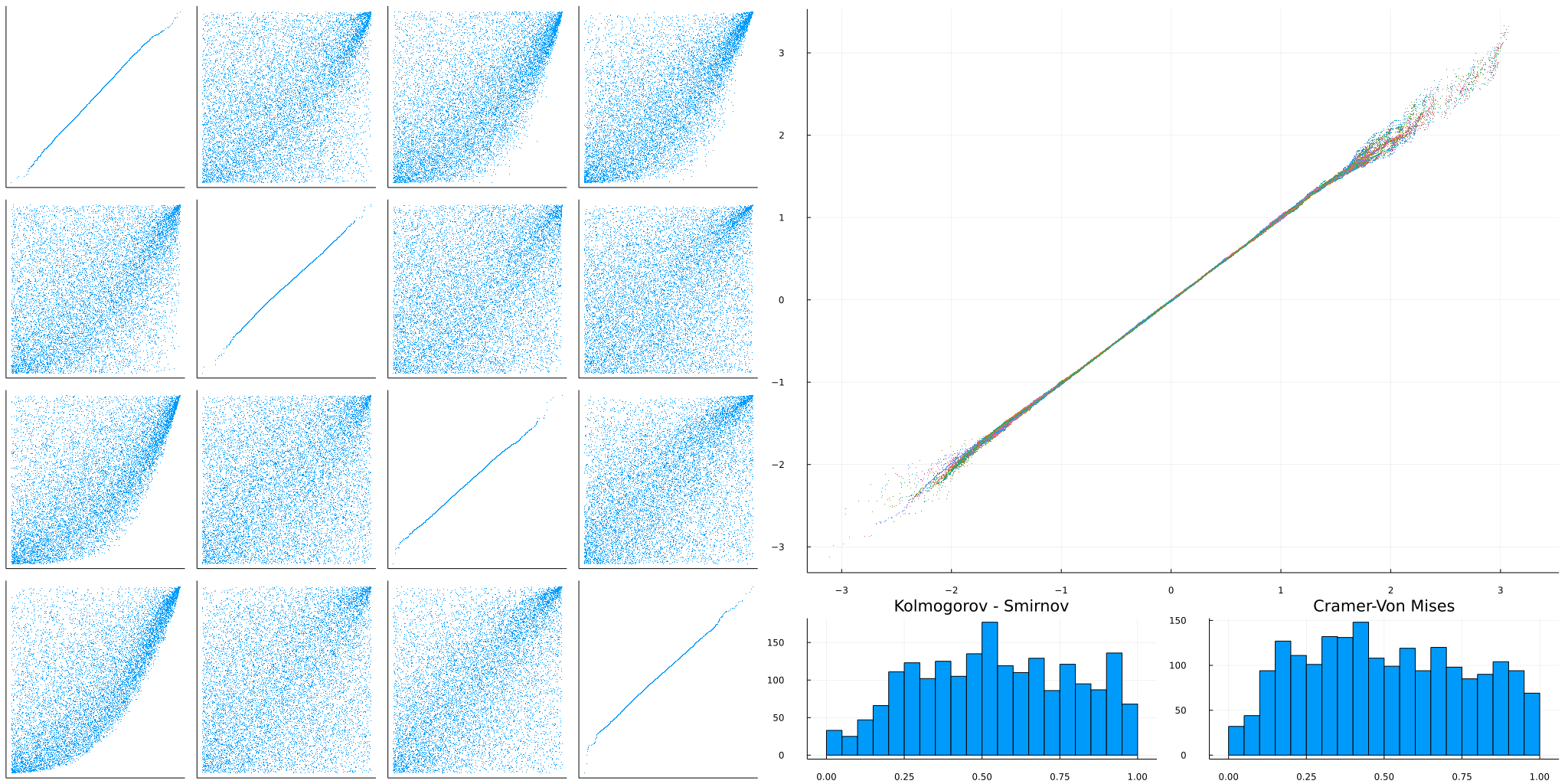}
  \caption{Results on Dataset~\ref{ggchd:data:WeirdExample} after thresholding the Thorin measure. On the left, a pairs plot representing the four dimensions: the bottom left triangle represents the pseudo-data from the original dataset, while the top triangle represent the pseudo-data from our estimator, and the middle line represents quantile-quantile plots of the marginals. Top right, quantile-quantile plots of $\langle \bm c_i, \hat{\bm X}\rangle$ on randomly chosen directions $\bm c_1,...,\bm c_{50}$. Bottom right, histograms of approximated Kolmogorov-Smirnov and Cramer-Von~Mises p values of our estimation.}
  \label{ggchd:fig:WeirdExample_threasold_summary}
\end{figure}

\subsection{The Loss-Alae example}

The next example we provide is on the \emph{Loss-Alae} dataset from Klugman \& Parsa~\cite{klugman1999}. This dataset features $1500$ observations of two variables: “Loss” and “Alae”.
We refer to ~\cite{klugman1999} for description of this (quite standard) dataset
in insurance valuations and copula estimations fields. This dataset was already considered in~\cite{laverny2021a} through the bruteforce minimisation of the $\mathcal L$ loss, which yielded very good results but required high-precision computations and gradient-free optimisation routines. We show here that our approximation also provide a very good Gamma deconvolution for this dataset, while working on standard precision and without the need for global optimisation routines. 

We ran the stochastic gradient descent on $\widetilde{\mathcal L}$ with $n=500$ potential atoms for the Thorin measure. Results are given in Figure~\ref{ggchd:fig:lossalae}.

\figuretb
  \centering
  \includegraphics[width=\linewidth]{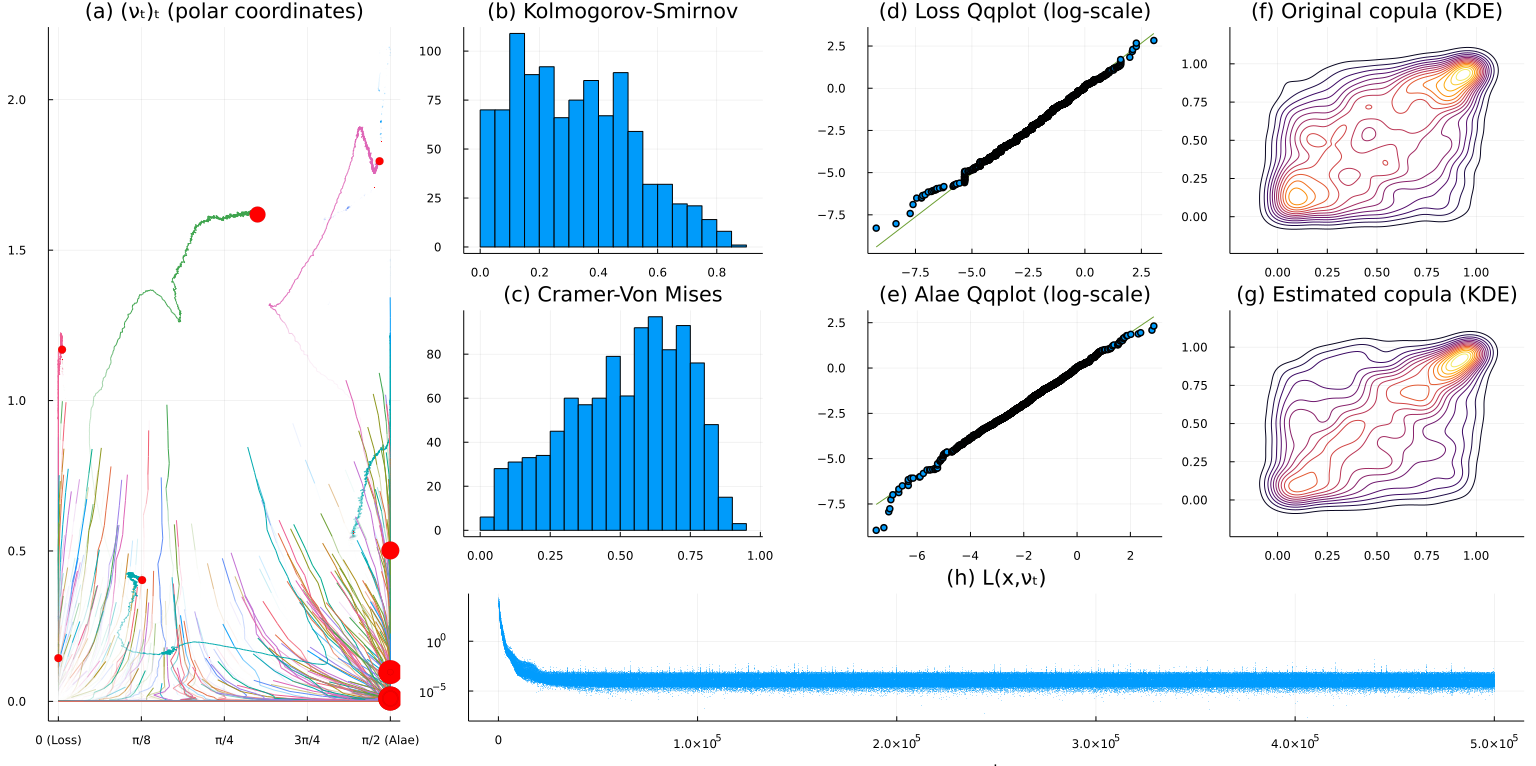}
  \caption{Results on the Loss-Alae dataset. (a) A trajectory of $\nu_t$ initialized on $\nu_0$ picked from $5000$ random atoms (Gaussian noise). The bivariate atoms are represented on polar coordinates (on the first quarter of the plane, Loss is on the left and Alae on the right), and the final weights are represented by the size of the red dots. (b, c) Kolmogorov-Smirnov and Cramer-Von Mises resampled p-values. (d, e) Quantile-quantile plots of the marginals. (f, g) Original pseud-observations and pseudo-observations taken from a simulation of the final estimator. (h) Loss $\widetilde{\mathcal L}(\bm x, \nu_t)$ on the $500000$ first steps.}
  \label{ggchd:fig:lossalae}
\end{figure}

As Figure~\ref{ggchd:fig:lossalae} (a) shows, most of the $n=500$ atoms annihilated and converged to zero. Only a few atoms are remaining, among which many have zero scales in one of the two dimension -- corresponding to red points at abscissa $0$ (Atoms are only on the Loss marginal) and $\frac{\pi}{2}$ (the Alae marginal), due to the polar coordinate system of this plot. Only eight red points are showing up on the graph, but the final measure has a few more very small atoms (with weights smaller than $10^{-5}$), which are not plotted. Remark the fact that most of the particles are either annihilated or converging towards one of the final remaining atoms. This is particularly true for scales that are zero in one of the two dimensions. On the other hand, scales that are joint between the two dimensions seems to be generated later in the process.
Our final estimator is given in Table~\ref{ggchd:tab:thorin_measure_LossAlae}.

\tabletb
  \scriptsize
	\caption{\label{ggchd:tab:thorin_measure_LossAlae}Thorin measure obtained on the Loss-Alae Dataset, after a thresholding at $10^{-5}$.}
	\centering
	\begin{tabu} to \linewidth {>{\raggedleft}l>{\raggedleft}X>{\raggedleft}X}
	\toprule
  Shapes    & Loss scales      & Alae scales\\
	\midrule
	\addlinespace[0.3em]
  $0.0001$  & $2.121$  & $0.0417$\\
  $0.0001$  & $0.1345$ & $0.7483$\\
  $0.0005$  & $2.0968$ & $0.042$ \\
  $0.0198$  & $0.2659$ & $4.9842$\\
  $0.0445$  & $2.2127$ & $0.0408$\\
  $0.0487$  & $2.3684$ & $3.2722$\\
  $0.1396$  & -        & $0.6481$\\
  $0.1739$  & $0.4553$ & $0.1914$\\
  $0.2544$  & -        & $0.1014$\\
  $0.7765$  & $0.1544$ & -       \\
  $0.8464$  & -        & $0.0092$\\
	\bottomrule
	\end{tabu}
\end{table}

The resulting estimator reproduces fairly well the dependence structure, see Figure~\ref{ggchd:fig:lossalae} (f, g), and the marginal quantile-quantile plots on Figure~\ref{ggchd:fig:lossalae} (d, e) are also quite good. The Kolmogorov-Smirnov and Cramer-Von Mises resampled $p$-values (see Appendix~\ref{ggchd:apx:KSCvM}) confirm the fact that we did fit correctly the dataset.
In terms of performance, obtaining this estimator took less than 5 seconds on a standard 2018 laptop.

\subsection{High dimensional multiplicative model}

As finding a high dimensional public dataset that exhibits features corresponding to the class is complicated\footnote{There is a typical case that is of important interest to us: the internal modeling in Non-life (re)-insurance modeling. See, e.g., \cite{arbenz2012,ferriero2011} for details on this case and on why a convolutional structure makes sense. However, such datasets are usually private.}, we propose to construct one. If we know that $\mathcal G_1$ is closed w.r.t. product of independent random variables~\cite{bondesson2015}, no such thing is known for the multivariate extensions. Therefore, we propose in Dataset~\ref{ggchd:data:BigExample} a multiplicative structure. 

\begin{dataset}[High dimensional multiplicative data]\label{ggchd:data:BigExample} For $d=2000$, let $G \sim\mathcal G_{1}(\delta_{1})$, $H \sim\mathcal G_{1}(2\delta_{\frac{1}{2}})$, and $Z_1,...,Z_d \sim \mathcal N(0,1)$ be all independent random variables, and let $\alpha_1,..,\alpha_d$ be fixed parameters in $[0,1]$ (uniformly chosen). Construct the random vector $\bm X = (X_1,...,X_d)$ as: $$\bm X = \left(G\,e^{Z_i}\,H^{1+2\alpha_i}\right)_{i \in 1,...,d}.$$
  We simulate a dataset $\bm x \in \mathbb R^{1500 \times 2000}$ of $N = 1500$ i.i.d. samples from $\bm X$.
\end{dataset}

The construction of Dataset~\ref{ggchd:data:BigExample} ensures that every marginal is inside the $\mathcal G_1$ class since the class is closed w.r.t. multiplication of random variables, contains log-Normals and powers of Gammas. However, we have no clue about whether the constructed distribution belongs to $\mathcal G_{2000}$ or not due to the dependency induced by the multiplicative risk-factors $G$ and $H$. Furthermore, there are only $1500$ observations, and hence the data is quite sparse. You can observe in Figure~\ref{ggchd:fig:BigExample_heatmap} Pearson, Spearman, and Kendall correlations heatmaps for this dataset, which all show that the dependency is quite strong (although in some dimensions clearly non-linear). 

\figuretb
  \centering
  \includegraphics[width=\linewidth]{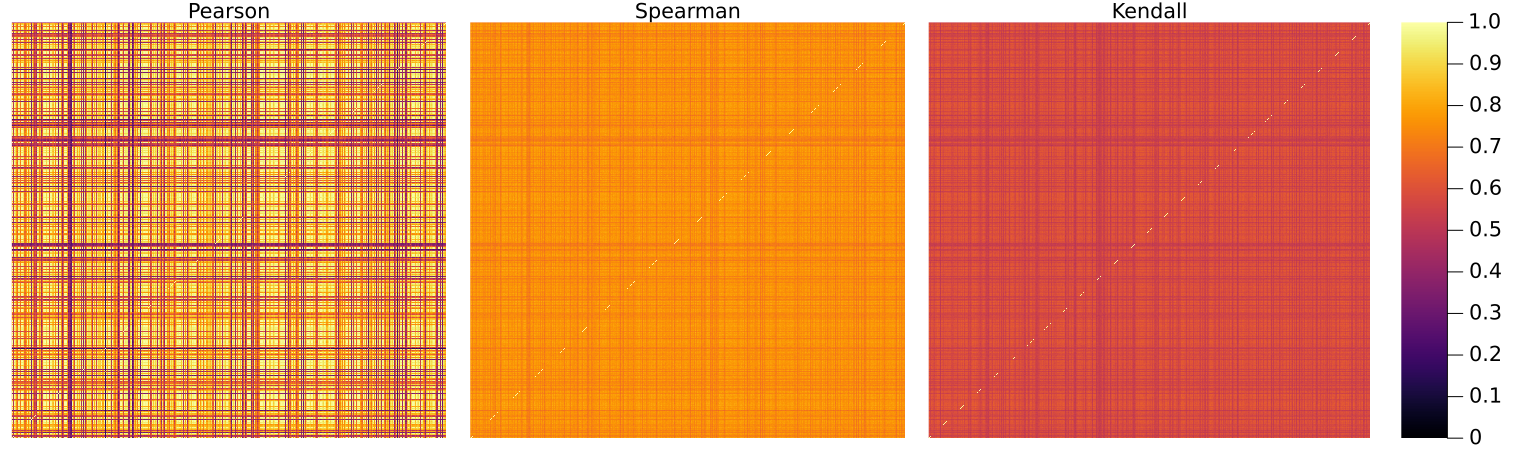}
  \caption{Heatmaps of pairwise correlation coefficients for Dataset~\ref{ggchd:data:BigExample}.}
  \label{ggchd:fig:BigExample_heatmap}
\end{figure}

We nevertheless launch the same gradient descent as before, with $n=500$ atoms only to see if we can obtain a relevant Thorin measure. The iteration cost is linear in $N,d$ and $n$, which allows dealing with these kinds of dimensions without issues on a standard laptop. 

Unfortunately, the empirical version of Kolmogorov-Smirnov and Cramer-Von Mises tests from Appendix~\ref{ggchd:apx:KSCvM} are not possible in these settings, as they would require a very big number of samples for their benchmarks (of some kind of importance sampling). Pairs plots are also not possible for the whole dataset (only for subsets of the dimensions), but we can still check the adequateness of our model through quantile-quantile plots on random projections, which are given in Figure~\ref{ggchd:fig:BigExample_summary}.

\figuretb
  \centering
  \includegraphics[width=\linewidth]{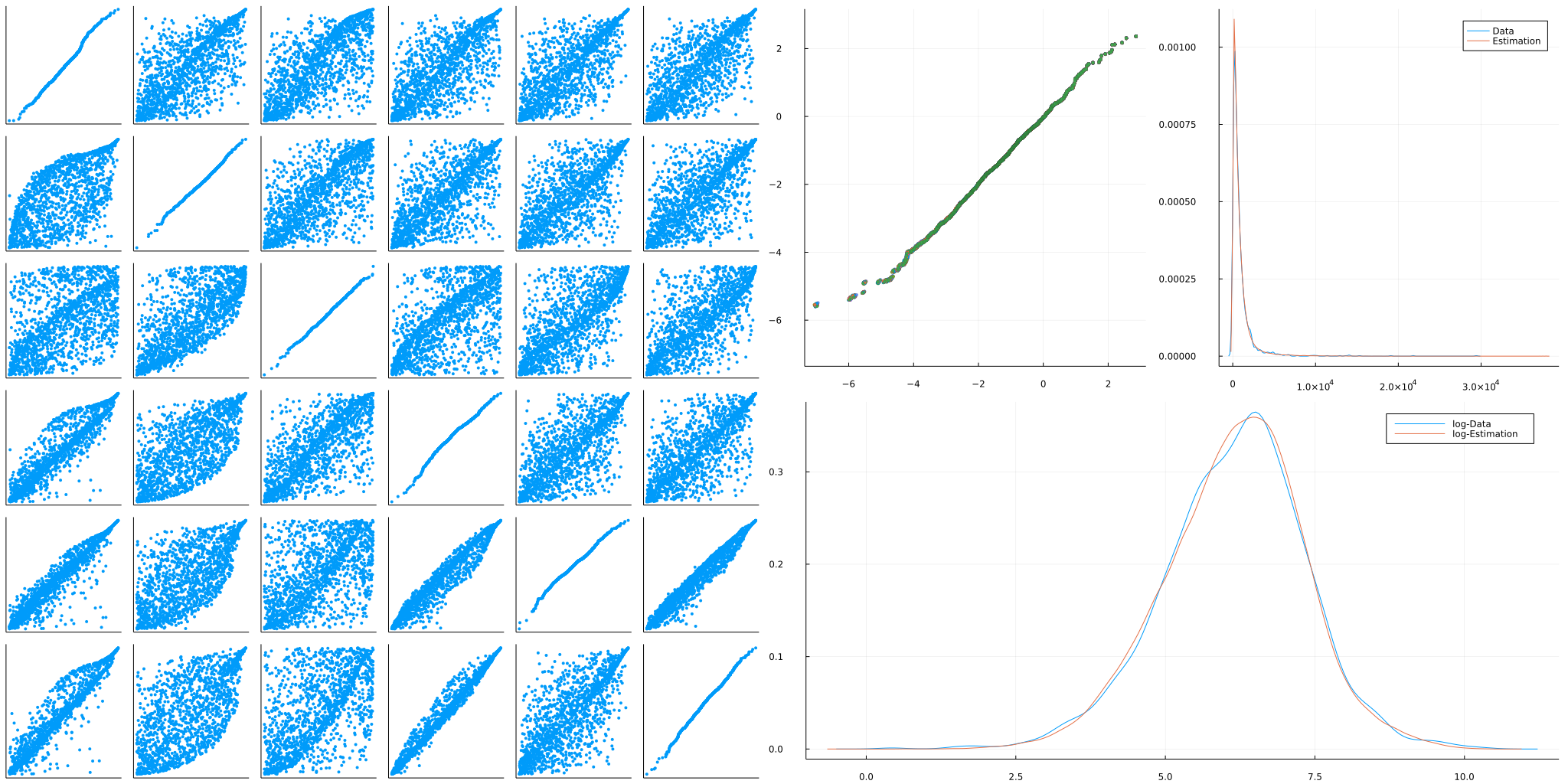}
  \caption{On the left, pairs plots of the 6 randomly chosen dimensions; On the top-right, a quantile-quantile plot of 50 random chosen projections, and the densities of $\langle \bm 1, \bm x\rangle$ and $\langle \bm 1, \hat{\bm X}\rangle$, showing the fit on the global sum, where $\bm x$ is the input data and $\hat{\bm X}$ is the estimated distribution. The lower-right density plots correspond to the densities of $\ln(\langle \bm 1, \bm x\rangle)$ and $\ln(\langle \bm 1, \hat{\bm X}\rangle)$. The densities of the univariate Gamma convolutions $\ln(\langle \bm 1, \hat{\bm X}\rangle)$ are computed by the Moschopoulos algorithm~\cite{moschopoulos1985}.}
  \label{ggchd:fig:BigExample_summary}
\end{figure}

The results on Dataset~\ref{ggchd:data:BigExample} were obtained with only $10^6$ iterations of our gradient descent, which means that only $10^6$ linear projections contributed one gradient each to the gradient descent. On the other hand, the constant $D(m,d)$ is here of order $10^{30}$, which means that we only explored a negligible portion of the necessary linear projections to recover the multivariate Thorin moments through Property~\ref{ggchd:prop:linearisation_moments}.

The projected quantile-quantile plot is almost perfect, and the results on the sum $\langle \bm 1, \bm x\rangle$ are correct. However, the details of the dependence structure are not perfectly fitted. If we let the algorithm run longer, it gets better. Finally, we wanted to say that only $160$ atoms are remaining (i.e., have weights bigger than $1e-16$), among which only $31$ have weights bigger than $1e-5$. The Thorin measure is very sparse, as much of the atoms are not shared by all dimensions and the matrix $\bm s$ of positions of the $160$ atoms has $40\%$ of its values that are zeros, and $20\%$ more are smaller than $1e-5$.

\section{Conclusion}\label{ggchd:sec:conclusion}
\ifinthesis\addtocontents{toc}{\protect\vskip-11pt}\fi

While the class of multivariate generalized convolutions is quite peculiar, it has a lot of interesting properties for many application fields. The first available estimator for distributions in this class, provided by~\cite{laverny2021a}, had good convergence properties, but lacked usability in higher dimensions.

We showed that fitting multivariate Gamma convolutions in high dimensions was nevertheless possible thanks to the convolutional properties of the class, that allow to bypass the dimensionality issues through random projections and obtain an efficient stochastic conic particle gradient descent. We showed first the consistency of the obtained global minimizer. Then, although our loss is non-convex, we showed the convergence of the associated gradient flow to a global minimizer. 

The behavior of the conic particle gradient descent is really promising, and the inclusion of a lasso penalty on the Thorin measure is quite easy to do in our settings. Although the theoretical results ensure the convergence of the gradient flow, our numerical section uses a stochastic approximation of the discretization of this flow, on which little is known. We nevertheless have good numerical results. 

However, the choice of the expansion size in each projection, represented by the hyperparameter $m$, is still something that should be considered. This parameter represents a level of \emph{details} that we try to learn from the data (but not the capacity of the resulting model to acknowledge these details, which would be represented by the number of atoms $n$ of the Thorin measure), and therefore should probably be learned through a cross-validation method. Inclusion of a lasso constraint from Remark~\ref{ggchd:rem:lasso} could be carried on via the same cross-validation rounds. }

\begin{appendix}
    {\graphicspath{{img/}} 
\section{Detailed gradients' computation}\label{ggchd:apx:nabla_AB}
\ifinthesis\addtocontents{toc}{\protect\vskip-11pt}\fi

In this appendix we describe more deeply what is the function $\bm B$, and how we perform computation of the gradients of $\widehat{\mathcal L}(\langle \bm c,\xdata\rangle, \cdot_{\langle \bm c\rangle})$. See the code of \texttt{ThorinDistributions.jl} for deeper implementation details. 

Recall that the function $\bm B$ is defined implicitly from the equation $\bm \mu = \bm B(\bm \tau)$, where $\bm\mu$ are the Taylor coefficients of a function $M$ and $\bm\tau$ are (up to scaling) Taylor coefficients of the function $K$ such that $M = \exp \circ K$. More precisely:
$$\sum_{\bm k \in \mathbb N^d} \frac{\mu_{\bm k}}{\bm k !} (\bm t - \bm 1)^{\bm k} = \exp\left\{\sum_{\bm k \in \mathbb N^d} \frac{\left(\lvert \bm k \rvert - 1\right)! \tau_{\bm k}}{\bm k !} (\bm t - \bm 1)^{\bm k}\right\}.$$

Up to the scaling by $\left(\lvert \bm k \rvert - 1\right)!$, the link between $\bm \mu$ and $\bm \tau$ is therefore directly given by multivariate Bell polynomials (see~\cite{bernardini2005,cvijovic2011,kim2019,melman2018}), from which we extract the combinatorial expressions for the function $\bm B$ given in Property~\ref{ggchd:prop:bijection_mu_kappa}.

\begin{property}[Bijection between moments and cumulants]\label{ggchd:prop:bijection_mu_kappa}
	For a given increasing index set $I \subseteq \mathbb N^d$, 
	  the sets $\left(\mu_{\bm i}\right)_{\bm i \in I}$ and $\left(\tau_{\bm i}\right)_{\bm i \in I}$ can be expressed from one another through:
	
	\begin{align*}
	\mu_{\bm i} &= B_{\bm i}(\bm\tau) =  e^{\tau_{\bm 0}}\sum\limits_{\pi \in \Pi(\bm i)} \prod\limits_{B \in \pi} \tau_{\bm i(B)}(\lvert \bm i(B) \rvert -1)!\\
	\tau_{\bm i} &= B^{-1}_{\bm i}(\bm \mu) = \frac{1}{(\lvert \bm i(B) \rvert -1)!}\sum\limits_{\pi \in \Pi(\bm i)} (\lvert \pi \rvert - 1)! (-1)^{\lvert \pi \rvert - 1} \left(\mu_{\bm 0}\right)^{-\lvert \pi \rvert} \prod\limits_{B \in \pi} \mu_{\bm i(B)}
	\end{align*}
	
	where $\Pi(\bm i) = \mathcal P(B(\bm i))$ is the set of all partitions of the multi-set $$B(\bm i) = \bigcup\limits_{j=1}^{d} \bigcup\limits_{p=1}^{i_j} \{j\},$$
	  while a partition $\pi$ in $\Pi(\bm i)$ is constituted of multi-sets $B$ composed of integers,
	  and $$\bm i(B) = \left(\#\{\text{occurences of }j\text{ in }B\}\right)_{j \in \{1,...,d\}}.$$
	
	\begin{proof}
		This is a direct application of the multivariate Faà Di Bruno's formula~\cite{constantine1996,mishkov2000,nardo2011}
	  to $M = \exp \circ K$ in one hand and $K = \ln \circ M$ on the other hand.
    The only details worth noting is that $\bm\tau$ are not exactly the derivatives of $K$ at $-\bm 1$, but rather $\tau_{\bm k} (\lvert \bm k \rvert -1)!$ are, according to Definition~\ref{ggchd:def:tau} and Property~\ref{ggchd:prop:thorin_moments}.
	\end{proof}

\end{property}

\begin{remark}[$\bm B$ and $\bm B^{-1}$ are almost polynomials] If we consider that $\mu_{\bm 0} = e^{\tau_{\bm 0}}$ is a fixed constant, the function $\bm B$ and its inverse are both polynomials in their remaining parameters.
\end{remark}

In the direct sense, i.e., to compute $\bm \mu = \bm B(\bm \tau)$,
  the good property of the exponential function of being its own derivative allows a recursive approach,
  adapted from~\cite{miatto2019}, described in Algorithm~\ref{ggchd:algo:Miatto2019}.

\algohtb
	\SetAlgoLined
	\KwIn{An increasing set $I$, Thorin moments $\bm \tau = \left(\tau_{\bm k}\right)_{\bm k \in I}$}
	\KwResult{Shifted moments $\bm \mu = \left(\mu_{\bm k}\right)_{\bm k \in I} = \bm B(\bm \tau)$}
	Set $\mu_{\bm 0} = \exp\left(\tau_{\bm 0}\right)$\\
	Set $\mu_{\bm k} = 0$ for all $\bm k \neq \bm 0$\\
	\ForEach{$\bm k:\;\bm 0 \neq \bm k \in I$}{ 
		Set $i_0$ as the index of the first element of $\bm k$ that is non-zero.\\
		Set $\bm p = \bm k$\\
		Set $p_{i_0} = p_{i_0} - 1$\\
		\ForEach{$\bm \ell:\; \bm \ell \le \bm p$}{
			$\mu_{\bm k} \pluseq \mu_{\bm \ell}  \tau_{\bm k - \bm \ell}  \binom{\bm p}{\bm \ell} (\bm k - \bm \ell)!$
		}
	}
	Return $\bm \mu$
	\caption{Recursive computation of the function $\bm B$. (see~\cite{miatto2019})}
  \label{ggchd:algo:Miatto2019}
\end{algorithm}

Algorithm~\ref{ggchd:algo:Miatto2019} is notable because it computes $\bm\mu$ from $\bm\tau$ in a time that is mostly quadratic in the number of coefficients $\lvert I \rvert$, while the raw usage of the Multivariate Faà Di Bruno expression given in Property~\ref{ggchd:prop:bijection_mu_kappa} has an exponential complexity. This increase in speed is due to the fact that the link function is the exponential function ($M(\bm t) = \exp\left\{K(\bm t)\right\}$), and is not possible for other link functions.  

We now expose the computation of the gradient of $\widehat{\mathcal L}(\langle\bm c,\xdata\rangle,\nu_{\langle \bm c\rangle})$, when $\nu = \sum_{i=1}^n p_i^2 \delta_{(q_{i,j}^2)_{j}}$, with respect to $\bm p, \bm q$. This square root parametrization allows using unconstrained optimization which is desired for numerical purposes. Recall that : 

$$\widehat{\mathcal L}(\langle\bm c,\xdata\rangle,\nu_{\langle \bm c\rangle}) := \lVert \hat{\bm\tau}(\langle\bm c,\xdata\rangle) - \bm \tau(\nu_{\langle \bm c\rangle}) \rVert_{\nabla(\langle\bm c,\xdata\rangle)}^2$$

Since only $\bm \tau\left(\nu_{\langle \bm c\rangle}\right)$ depends on the parameters, we have that : 

\begin{equation}\label{ggchd:eq:gradient_expression}
  \nabla_{\widehat{\mathcal L}}(\bm p,\bm q) = 2\left(\hat{\bm\tau}\left(\langle\bm c,\xdata\rangle\right) - \bm \tau\left(\sum_{i=1}^n \alpha_i \delta_{\langle\bm c, \bm s \rangle}\right)\right)'\nabla\left(\langle\bm c,\xdata\rangle\right)'\nabla\left(\langle\bm c,\xdata\rangle\right) \nabla_{\tau\left(\sum_{i=1}^n \alpha_i \delta_{\langle\bm c, \bm s \rangle}\right)}\left(\bm p, \bm q\right)
\end{equation}

Furthermore, recall that for an atomic measure $\nu$, 
$$\tau_{\bm 0}(\nu) = -\sum_{i=1} \alpha_i \ln\left(1+\lvert\bm s_i\rvert\right)  \text{ and }\tau_{\bm k}(\nu) = \sum\limits_{i=1}^n \alpha_i \left(\frac{\bm s_i}{1 + \lvert\bm s_i\rvert}\right)^{\bm k} \text{ for } \bm k \in \mathbb N_*^d.$$

Hence, $\nabla_{\tau\left(\sum_{i=1}^n \alpha_i \delta_{\langle\bm c, \bm s \rangle}\right)}(\bm p, \bm q)$ is composed of :

\begin{align*}
  \frac{\partial \tau_{0}}{\partial p_i} &= -2p_i \ln\left(1 + \langle \bm c, \bm q_i^2\rangle\right)\\
  \frac{\partial \tau_{0}}{\partial q_{i,j}} &= 2p_i \left(\frac{\langle \bm c, \bm q_i^2\rangle}{1 + \langle \bm c, \bm q_i^2\rangle}\right)^k\\
  \frac{\partial \tau_{k}}{\partial p_i} &= \frac{2p_i^2 q_{i,j} c_j}{1 + \langle \bm c, \bm q_i^2\rangle}\\
  \frac{\partial \tau_{k}}{\partial q_{i,j}} &= k \left(\frac{\langle \bm c, \bm q_i^2\rangle}{1 + \langle \bm c, \bm q_i^2\rangle}\right)^{k-1} \frac{2p_i^2 q_{i,j} c_j}{\left(1 + \langle \bm c, \bm q_i^2\rangle\right)^2}\\ 
\end{align*}

where in the last line, the factor $p_i^2$ can be removed to ensure the compatibility with the metrisation of Section~\ref{ggchd:sec:gradient_flows} (independence of spatial and vertical components of the gradient). 

Moreover, $\nabla_{\bm A\bm B}(\bm \tau(\langle\bm c,\xdata\rangle)) = \bm A\nabla_{\bm B}(\bm \tau(\langle\bm c,\xdata\rangle))$, and we compute $\bm J = \nabla_{\bm B}(\bm \tau(\langle\bm c,\xdata\rangle))$ the second part through the Algorithm~\ref{ggchd:algo:derivative}, which is constructed through a simple chain rule run through Algorithm~\ref{ggchd:algo:Miatto2019}. 

\algohtb
	\SetAlgoLined
	\KwIn{An increasing set $I$, Thorin moments $\bm \tau = \left(\tau_{\bm k}\right)_{\bm k \in I}$}
	\KwResult{The Jacobian $\bm J \in \mathbb R^{\lvert I\rvert \times \lvert I \rvert}$ of $\bm B$ taken at $\bm\tau$}
  Set $J_{\bm 0,\bm 0} = \exp\left(\tau_{\bm 0}\right)$.\\
	\ForEach{$\bm k:\;\bm 0 \neq \bm k \in I$}{ 
		Set $i_0$ as the index of the first element of $\bm k$ that is non-zero.\\
		Set $\bm p = \bm k$\\
		Set $p_{i_0} = p_{i_0} - 1$\\
		\ForEach{$\bm l:\; \bm l \le \bm p$}{
      $c = \binom{\bm p}{\bm l} (\bm k - \bm l)!$\\
      $J_{\bm k, \cdot} \pluseq cJ_{\bm l, \cdot} \tau_{\bm k - \bm l}$\\
      $J_{\bm k, \bm k - \bm l} \pluseq c\mu_{\bm l}$\\
			$\mu_{\bm k} \pluseq c \mu_{\bm l}  \tau_{\bm k - \bm l}$
		}
	}
	Return $\bm J$
	\caption{Recursive computation of the function $\nabla_{\bm B}$.}
  \label{ggchd:algo:derivative}
\end{algorithm}

Note that Algorithm~\ref{ggchd:algo:derivative} also computes $\mu$, which is necessary due to the recursive nature of the chain rule, and therefore computations of Algorithms~\ref{ggchd:algo:Miatto2019} and~\ref{ggchd:algo:derivative} can (and should) be pooled.

Last but not least, the computation of $\bm \tau(\langle \bm c, \xdata \rangle)$ is worth mentioning. According to Definition~\ref{ggchd:def:mc_estimators}, we have : 
    $$\hat{\bm\tau}(\xdata) = \bm B^{-1}\left(\left(\frac{1}{N}\sum_{i=1}^N \xdatai^{\bm k} e^{-\lvert \xdatai \rvert}\right)_{\bm k \in I}\right).$$
and we thus require the computation of $\bm B^{-1}$. If $\bm B^{-1}$ is a very costly function to compute in any dimension, as Property~\ref{ggchd:prop:bijection_mu_kappa} might show, it however drastically simplifies in dimension one, and its expression can be translated into the efficient Algorithm~\ref{ggchd:algo:invB_dim1}.

\algohtb
	\SetAlgoLined
	\KwIn{$m \in \mathbb N$ and a univariate dataset $\xdata \in \mathbb R^{N}$}
	\KwResult{$\bm\tau = \left(\tau_k(\xdata)\right)_{k \le m}$}
  Compute first $\bm \mu = \left(\frac{1}{N}\sum_{i=1}^N x_i^{k} e^{-x_i}\right)_{k \le m}$ any way you want (e.g., via recursion on $k$.) \\
  Let $\bm\eta = \left(\frac{\mu_k}{k!\mu_0}\right)_{k\le m}$\\
  Set $\tau_{0} = \ln\left(\mu_0\right)$\\
	\ForEach{$k \in \{1,...,m\}$}{ 
    $\tau_k = (k-1) * \eta_k$\\
		\ForEach{$j \in \{1,...,k-1\}$}{
      $\tau_k \moinseq \tau_j * \eta_{k-j}$
		}
	}
	Return $\bm \tau$
	\caption{Computation of Thorin moments of a univariate dataset.}
  \label{ggchd:algo:invB_dim1}
\end{algorithm}
  
Note that although Algorithm~\ref{ggchd:algo:invB_dim1} has a quadratic complexity in the number of coefficients, unlike Algorithm~\ref{ggchd:algo:Miatto2019}, this will not hold in higher dimensions $d > 1$ (the log function in $K(t) = \ln\left\{M(t)\right\}$ is not as nice as the $\exp$ function), which is one of the reasons we use univariate projections. Our renormalization of the derivatives of $K(t)$ also simplifies the inner loop of Algorithm~\ref{ggchd:algo:invB_dim1} by removing a multiplication by a binomial coefficient, reducing again the computational cost. 

To obtain the stochastic gradient of the global loss, the only remaining computations are matrix-vector multiplications (that must be done in the right order for performance), according to Equation \eqref{ggchd:eq:gradient_expression}.

\section{Non-distribution free tests}\label{ggchd:apx:KSCvM}
\ifinthesis\addtocontents{toc}{\protect\vskip-11pt}\fi

Our goal here is to produce a test from the Kolmogorov-Smirnov (resp Cramer-Von~Mises) statistic in high dimension. 
The major problem is that these tests are not distribution free and the distribution of the statistics will depend on the distribution under the null hypothesis. However, through resampling, we are still able to construct a test that is well-calibrated. 

To do that, for a given random vector $\bm X$ whose distribution is known, 
we compute benchmark vectors $\mathrm{KS}_{\bm X}$ and $\mathrm{CvM}_{\bm X}$ that approximate the distributions of the statistics. This is detailed in Algorithm~\ref{ggchd:algo:benchmark}.  

First, recall that for a dataset $\bm D \in \mathbb R^{N\times d}$, we define the empirical distribution function as $$F_{\bm D} : \bm x \mapsto \frac{1}{N}\sum_{i=1}^N \mathbf 1_{D_{i,j} \le x_j \forall j \in 1,...,d}.$$

\algohtb
	\SetAlgoLined
	\KwIn{The distribution of a $d$-variate random vector $\bm X$, an integer $M$, a sample size $N$.}
	\KwResult{Samples $\mathrm{KS}_{\bm X}$ and $\mathrm{CvM}_{\bm X}$ from the Kolmogorov-Smirnov and Cramer-Von~Mises statistics distributions associated to $\bm X$.}
	\ForEach{$i \in 1,...,M$}{ 
		Sample a pair of datasets $\bm D_1,\bm D_2$ in $\mathbb R^{N \times d}$ from $\bm X$. \\
    Let $\mathrm{KS}_{\bm X,i} = \sup\limits_{x \in \bm D_1} \lvert F_{\bm D_1}(x) - F_{\bm D_2}(x) \rvert.$\\
    Let $\mathrm{CvM}_{\bm X,i} = \frac{1}{N}\sum_{x \in \bm D_1} (F_{\bm D_1}(x) - F_{\bm D_2}(x))^2 dx.$\\
	}
	Return $\mathrm{KS}_{\bm X}$ and $\mathrm{CvM}_{\bm X}$
	\caption{Kolmogorov-Smirnov and Cramer-Von~Mises distribution approximation of $\bm X$.}
  \label{ggchd:algo:benchmark}
\end{algorithm}

Then, for a proposal distribution $\hat{\bm X}$, we compute the statistics and the approximated p values through Algorithm~\ref{ggchd:algo:ks_cm_stats}.

\algohtb
	\SetAlgoLined
	\KwIn{The distributions of two $d$-variate random vectors $\bm X$ and $\hat{\bm X}$, a sample size $N$.}
	\KwResult{A sample of the KS and CvM distances between the two distributions.}
  Compute $\mathrm{KS}_{\bm X}$ and $\mathrm{CvM}_{\bm X}$ via Algorithm~\ref{ggchd:algo:benchmark}.\\ 
  Sample a pair of datasets $\bm D_1,\bm D_2$ in $\mathbb R^{N \times d}$ from $\bm X$ and $\hat{\bm X}$ respectively. \\
  Let $\mathrm{KS} = \sup\limits_{x \in \bm D_1} \lvert F_{\bm D_1}(x) - F_{\bm D_2}(x) \rvert.$\\
  Let $\mathrm{CvM} = \frac{1}{N}\sum_{x \in \bm D_1} (F_{\bm D_1}(x) - F_{\bm D_2}(x))^2 dx.$\\
	Return $\frac{1}{M}\sum\limits_{i=1}^M \mathrm{KS} \le \mathrm{KS}_{\bm X,i}$ and $\frac{1}{M}\sum\limits_{i=1}^M \mathrm{CvM} \le \mathrm{CvM}_{\bm X,i}$
	\caption{Approximated Kolmogorov-Smirnov and Cramer-Von~Mises p values.}
  \label{ggchd:algo:ks_cm_stats}
\end{algorithm}

These two algorithm give us a way to check if our estimators indeed are close to the distribution we simulated the dataset from. However, they are not practical as the computation of $\mathrm{KS}_{\bm X}$ and $\mathrm{CvM}_{\bm X}$ is quadratic in the number of samples in the dataset $N$ (although linear in the dimension $d$ and the number of resamples $M$), which can be dramatic in terms of runtime. 

As an example, we give a few graphs. Let $\bm X$ be a multivariate random vector with independent log-Normal(0,1) marginals and $\hat{\bm X} = \bm X + 2$. A simple run of the procedures gives results depicted in Figure~\ref{ggchd:fig:ks_test}.

\figuretb
  \centering
  \begin{subfigure}{.5\textwidth}
    \centering
    \includegraphics[width=\linewidth]{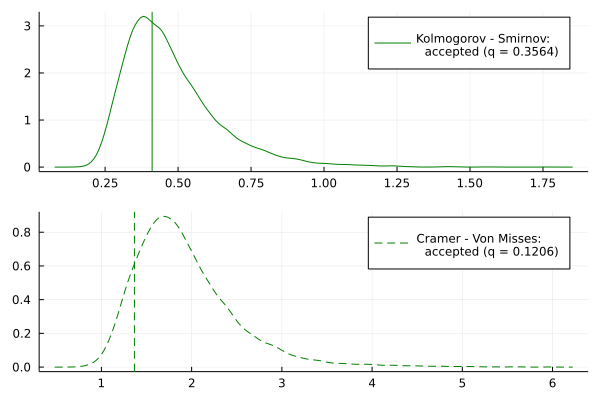}
    \caption{$\hat{\bm X} = \bm X$}
    \label{ggchd:fig:ks_test_accepted}
  \end{subfigure}%
  \begin{subfigure}{.5\textwidth}
    \centering
    \includegraphics[width=\linewidth]{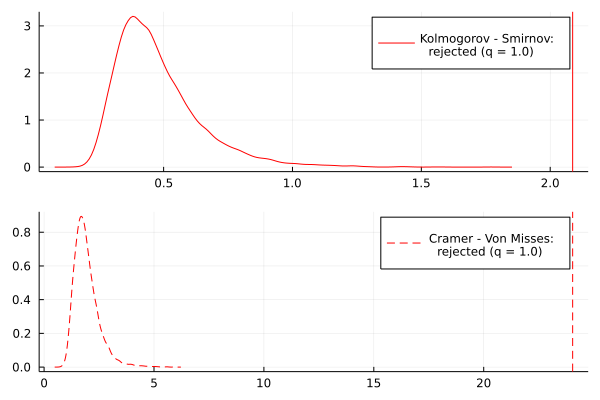}
    \caption{$\hat{\bm X} = \bm X + 2$}
    \label{ggchd:fig:ks_test_rejected}
  \end{subfigure}%
  \caption{Density of the (approximated) distributions of $\mathrm{KS}_{\bm X}$ (on the top) and $\mathrm{CvM}_{\bm X}$ (on the bottom). The vertical lines show values obtained by the estimator $\hat{\bm X}$. On the left, we are under the null and $\hat{\bm X}$ has the same distribution as $\bm X$, while not on the right.}
  \label{ggchd:fig:ks_test}
\end{figure}

Furthermore, if we resample the test, we obtain a histogram of p values, which are supposed to be approximately uniformly distributed under the null. Histograms of these p values are given in Figures~\ref{ggchd:fig:ks_test_histogram_good} and \ref{ggchd:fig:ks_test_histogram_bad} for $\bm X$ and $\hat{\bm X}$ respectively. 

\figuretb
  \centering
  \begin{subfigure}{.5\textwidth}
    \centering
    \includegraphics[width=\linewidth]{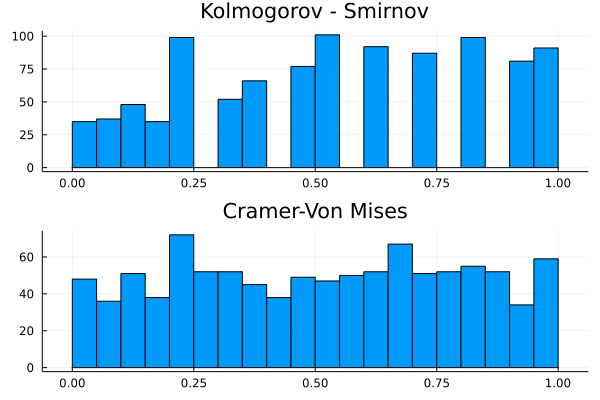}
    \caption{$\hat{\bm X} = \bm X$}
    \label{ggchd:fig:ks_test_histogram_good}
  \end{subfigure}%
  \begin{subfigure}{.5\textwidth}
    \centering
    \includegraphics[width=\linewidth]{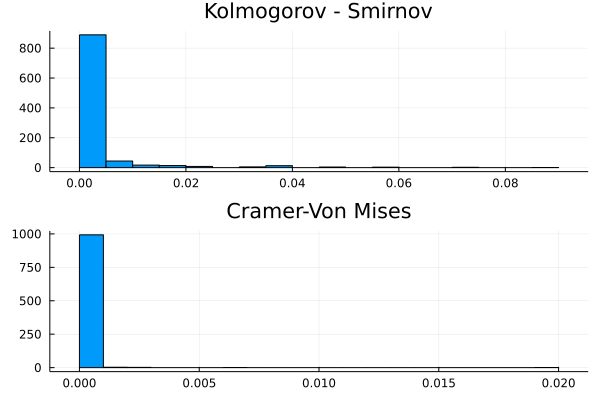}
    \caption{$\hat{\bm X} = \bm X + 2$}
    \label{ggchd:fig:ks_test_histogram_bad}
  \end{subfigure}%
  \caption{Histogram of resampled approximated p values of our test. On the left, under the null, but not on the right. }
  \label{ggchd:fig:ks_test_histogram}
\end{figure}

We see on Figure~\ref{ggchd:fig:ks_test_histogram_bad} that the p values are concentrated near zero, which forces us to reject the null hypothesis. Moreover, the spreading of p values under the null hypothesis that we observe in ~\ref{ggchd:fig:ks_test_histogram_good} is a good result for our test. 

\section{Non-Gamma extensions}
\ifinthesis\addtocontents{toc}{\protect\vskip-11pt}\fi

Our algorithm has a potential for non-Gamma extensions, which we detail in this appendix.

Suppose we know $K$, the cumulant generating function of a positive and infinitely divisible random variable $U$, that might \emph{not} be in the GGC class. 
Recall that since $U$ is positive, $K$ is at least analytic around $-1$. Consider now the shape-scale multivariate family based on $U$, that is the family of random vectors defined by their cumulant generating functions: 

$$K_{\nu}(\bm t) = \int K(\langle \bm s, \bm t\rangle) \nu(d\bm s), \text{ for }\nu\in\mathcal M_+(\mathbb R_+^d)$$

Conditions such that $K_{\nu}$ is a bona fide $d$-variate cumulant generating function are yet unclear outside the Gamma case. However, when $\nu$ is finitely atomic, this is always the case. Taking $U \sim \mathrm{Exp}(1)$ gives rise to the family of multivariate generalized Gamma convolutions. Other choices of $U$ might be possible and interesting. 

\begin{remark}[Link with infinite divisibility] When $0\le\alpha\le 1$, the cumulant generating function $K_{\alpha\delta_1}$ is the cumulant generating function of the $\alpha$-piece of $U$, in the divisibility sense, which we denote $U_{(\alpha)}$: If $U_{(\alpha)}$ and $U_{(1-\alpha)}$ are independent, their sum is distributed as $U$. 
\end{remark}

If the function $K$ is analytic around $-\bm 1$, then our work can be extended to this new convolutional class. However, for practical implementation, it is useful that the derivatives of $K$ are simple functions, in the Gamma case these a simple rational functions.
This reasoning allows to extend our algorithm and convergence results to convolution classes based on a \emph{non-Gamma infinite divisible distribution} for $U$. 

A possibility with non-atomic Thorin measure would be to use the positive $\alpha$-stable distribution, with parameter with $0 < \alpha < 1$, given by the cumulant generating function $$K(t) = -(-t)^{\alpha}.$$ According to~\cite[Ex 3.2.1]{bondesson1992}, the corresponding Thorin measure has density $$u(t) = \frac{\alpha \sin(\alpha\pi)}{\pi} t^{\alpha-1}.$$
Hence, a convolutional class based on this model would allow easily non-atomic Thorin measures and still simulate them. 

\section{An extension of Bélisle, Massé \& Ranford's result.}\label{ggchd:apx:projection}
\ifinthesis\addtocontents{toc}{\protect\vskip-11pt}\fi

This appendix proposes a new result, based on~\cite{renyi1952,heppes1956}, about linear projections and Grassmannian cubatures that was not finally used in our arguments. Keeping notations from Section~\ref{ggchd:sec:random_projs}, Rényi and Heppes provide the following theorem:  

\begin{theorem}[Rényi~\cite{renyi1952} - Heppes~\cite{heppes1956}]\label{ggchd:thm:RenyiHeppes} Let $\nu$ be a discrete probability measure on $\mathbb R^d$ with a support made with exactly $k$ distinct atoms. Assume that $E_1,...,E_{k+1}$ are subspaces of $\mathbb R^d$ of respective dimensions $d_1,...,d_{k+1}$ such that no couple of them is contained in a hyperplane. Then for any measure $\mu$ in $\mathbb R^d$, $\nu = \mu$ if and only if $E_i \in \mathcal E_{d_i}(\nu,\mu)$ for every $1 \le i \le k+1$.
\end{theorem}

This means that $k+1$ distinct projections onto $d-1$ dimensions each are enough to find a $d$-dimensional probability measure with $k$ atoms. More detailed exposition can be found in e.g., Chafaï~\cite{chafai2021} and references therein. 
However, we were mostly interested in projection onto univariate subspaces, for computational reasons. Therefore, we developed a more quantitative statement based on a certain condition on subspaces, that extends to any projection dimensions and non-normalized measures a simple results from Bélisle, Massé \& Ranford~\cite[Sec. 6]{belisle1997}.

\begin{theorem}[Distance bound between two measures.]\label{ggchd:thm:distnce_bound} Let $\nu_1$ and $\nu_2$ be two positive discrete measures on $\mathbb R^d$, let $k \in 1,...,d-1$ and let $E_1,...,E_n \in \mathcal E_{k}\left(\nu_1,\nu_2\right)$ ($n$ can be greater than $d$).
  
If $\mathcal H_k : \;\forall\, J \subseteq \left\{1,...,n\right\} s.t. \lvert J \rvert = d+1-k, \bigcap_{j \in J}E_j^{\perp} = \{\bm 0\}$, then
  $$\sup\limits_{\bm x \in \mathbb R^d} \lvert \nu_1(\bm x) - \nu_2(\bm x)\rvert \le \frac{d-k}{n} \max\{\lvert\nu_2\rvert,\lvert\nu_1\rvert\}.$$
\begin{proof}

  The case $k=d-1$ was already proved, but only for probability measures, in Bélisle, Massé \& Ranford~\cite[Sec. 6]{belisle1997} as an extension of Theorem~\ref{ggchd:thm:RenyiHeppes}. In this case, the condition $\mathcal H_k$ reduces to the fact that the $E_i$ are distinct elements of $\mathcal E_{d-1}$.  However, their proof does not hold for $k \neq d-1$, nor for general measures. We propose an extension here.

  Fix $\bm x \in \mathbb R^d$, and denote $c = \nu_1(\bm x) - \nu_2(\bm x)$. Then for each $E_i \in E_1,...,E_n$, define $A_i(\bm x)$ by : 
  \begin{align*}
    A_i(\bm x) &= \pi_{E_i}^{-1}\left(\pi_{E_i}(\bm x)\right) \setminus \{\bm x\}\\
        &= \bm x + (\mathrm{Ker}(\pi_{E_i}) \setminus \{\bm 0\}) \\
        &= \bm x + E_i^{\perp} \setminus \{\bm x\}.
  \end{align*}

  Since $E_i \in \mathcal E_{k}\left(\nu_1,\nu_2\right)$, $\nu_1\left(\bm x + E_i^{\perp}\right)=\nu_2\left(\bm x + E_i^{\perp}\right)$ where $E_i^{\perp}$ is of dimension $d-k$. 
  Hence: 
  \begin{align*}
    \nu_1(A_i(\bm x)) &= \nu_1\left(\bm x + E_i^{\perp}\right) - \nu_1(\bm x)\\
    &= \nu_2\left(\bm x + E_i^{\perp}\right) - \nu_2(\bm x) -c\\
    &= \nu_2(A_i(\bm x)) - c
  \end{align*}

  Therefore, $nc = \sum\limits_{i=1}^n \nu_2(A_i(\bm x)) - \sum\limits_{i=1}^n \nu_1(A_i(\bm x)) \le \sum\limits_{i=1}^n \nu_2(A_i(\bm x))$ since $\nu_1$ is a positive measure. 

  Now, since $\mathcal H_k$ is true, any $\bm y\in\mathbb R^d \setminus \{\bm 0\}$ can appear in at most $d-k$ different $E_i^{\perp}$. Indeed, all $E_i^{\perp}$ are of dimension $d-k$, and if a point other than $\{\bm 0\}$ (which they already share) appear in more than $d-k$ of them, they would have a common dimension, which is not possible.

  Hence, every atom of the measure $\nu_2$ can appear in at most $d-k$ different $A_i$. Therefore, denoting $S$ the support of $\nu_2$, 
  \begin{align*}
    nc &\le \sum\limits_{i=1}^n \nu_2(A_i(\bm x)) \\
    &=\sum\limits_{i=1}^n \sum\limits_{\bm y \in S\cap A_i} \nu_2(\bm y) \\
    &= \sum\limits_{\bm y \in S} \nu_2(\bm y) \sum\limits_{i=1}^n \mathbb 1_{\bm y \in A_i}\\
    &\le \sum\limits_{\bm y \in S} \nu_2(\bm y) (d-k) \\
    &= (d-k) \lvert\nu_2\rvert
  \end{align*}

  And hence $c \le \frac{d-k}{n} \lvert\nu_2\rvert$. By symmetry of the argument, $\lvert c \rvert \le \frac{d-k}{n} \max\{\lvert\nu_2\rvert,\lvert\nu_1\rvert\}$.
\end{proof}
\end{theorem}

\begin{remark}\label{ggchd:rem:hypothesis_H} The condition $\mathcal H_k$ from Theorem~\ref{ggchd:thm:distnce_bound} has two nice special cases: 
  \begin{itemize}
    \item $\mathcal H_{d-1}\;\iff$ the $E_i$ are distinct elements of $\mathcal E_{d-1}$, as~\cite[Sec. 6]{belisle1997} noted.
    \item $\mathcal H_{1}\;\iff$ any subset of size $d$ from the set of vectors $\{\bm c_1,...,\bm c_n\}$ forms a basis of $\mathbb R^d$.  
  \end{itemize}
\end{remark}

}
\end{appendix}

\begin{funding}
The author was supported by SCOR SE through a CIFRE grant. 
\end{funding}



\bibliographystyle{imsart-number} 
\bibliography{rfm}       

\begin{thebibliography}{47}

\bibitem{ambrosio2005}
\begin{bbook}[author]
\bauthor{\bsnm{Ambrosio},~\bfnm{Luigi}\binits{L.}},
  \bauthor{\bsnm{Gigli},~\bfnm{Nicola}\binits{N.}} \AND
  \bauthor{\bsnm{Savar{\'e}},~\bfnm{Giuseppe}\binits{G.}}
(\byear{2005}).
\btitle{Gradient flows: in metric spaces and in the space of probability
  measures}.
\bpublisher{Springer Science \& Business Media}.
\end{bbook}
\endbibitem

\bibitem{arbenz2012}
\begin{barticle}[author]
\bauthor{\bsnm{Arbenz},~\bfnm{Philipp}\binits{P.}},
  \bauthor{\bsnm{Hummel},~\bfnm{Christoph}\binits{C.}} \AND
  \bauthor{\bsnm{Mainik},~\bfnm{Georg}\binits{G.}}
(\byear{2012}).
\btitle{Copula based hierarchical risk aggregation through sample reordering}.
\bjournal{Insurance: Mathematics and Economics}
\bvolume{51}
\bpages{122--133}.
\end{barticle}
\endbibitem

\bibitem{belisle1997}
\begin{barticle}[author]
\bauthor{\bsnm{B{\'e}lisle},~\bfnm{Claude}\binits{C.}},
  \bauthor{\bsnm{Mass{\'e}},~\bfnm{Jean-Claude}\binits{J.-C.}} \AND
  \bauthor{\bsnm{Ransford},~\bfnm{Thomas}\binits{T.}}
(\byear{1997}).
\btitle{When is a probability measure determined by infinitely many
  projections?}
\bjournal{The Annals of Probability}
\bpages{767--786}.
\end{barticle}
\endbibitem

\bibitem{bernardini2005}
\begin{barticle}[author]
\bauthor{\bsnm{Bernardini},~\bfnm{A.}\binits{A.}},
  \bauthor{\bsnm{Natalini},~\bfnm{P.}\binits{P.}} \AND
  \bauthor{\bsnm{Ricci},~\bfnm{P.~E.}\binits{P.~E.}}
(\byear{2005}).
\btitle{Multidimensional Bell Polynomials of Higher Order}.
\bjournal{Computers \& Mathematics with Applications}
\bvolume{50}
\bpages{1697--1708}.
\end{barticle}
\endbibitem

\bibitem{bondesson1992}
\begin{bbook}[author]
\bauthor{\bsnm{Bondesson},~\bfnm{Lennart}\binits{L.}}
(\byear{1992}).
\btitle{Generalized {{Gamma Convolutions}} and {{Related Classes}} of
  {{Distributions}} and {{Densities}}}.
\bseries{Lecture {{Notes}} in {{Statistics}}}
\bvolume{76}.
\bpublisher{{Springer New York}}, \baddress{{New York, NY}}.
\end{bbook}
\endbibitem

\bibitem{bondesson2009}
\begin{barticle}[author]
\bauthor{\bsnm{Bondesson},~\bfnm{Lennart}\binits{L.}}
(\byear{2009}).
\btitle{On Univariate and Bivariate Generalized Gamma Convolutions}.
\bjournal{Journal of Statistical Planning and Inference}
\bvolume{139}
\bpages{3759--3765}.
\end{barticle}
\endbibitem

\bibitem{bondesson2015}
\begin{barticle}[author]
\bauthor{\bsnm{Bondesson},~\bfnm{Lennart}\binits{L.}}
(\byear{2015}).
\btitle{A {{Class}} of {{Probability Distributions}} That Is {{Closed}} with
  {{Respect}} to {{Addition}} as {{Well}} as {{Multiplication}} of
  {{Independent Random Variables}}}.
\bjournal{Journal of Theoretical Probability}
\bvolume{28}
\bpages{1063--1081}.
\end{barticle}
\endbibitem

\bibitem{bondesson2018}
\begin{barticle}[author]
\bauthor{\bsnm{Bondesson},~\bfnm{Lennart}\binits{L.}} \AND
  \bauthor{\bsnm{Simon},~\bfnm{Thomas}\binits{T.}}
(\byear{2018}).
\btitle{Stieltjes Functions of Finite Order and Hyperbolic Monotonicity}.
\bjournal{Transactions of the American Mathematical Society}
\bvolume{370}
\bpages{4201--4222}.
\end{barticle}
\endbibitem

\bibitem{breger2017}
\begin{bincollection}[author]
\bauthor{\bsnm{Breger},~\bfnm{Anna}\binits{A.}},
  \bauthor{\bsnm{Ehler},~\bfnm{Martin}\binits{M.}},
  \bauthor{\bsnm{Gr{\"a}f},~\bfnm{Manuel}\binits{M.}} \AND
  \bauthor{\bsnm{Peter},~\bfnm{Thomas}\binits{T.}}
(\byear{2017}).
\btitle{Cubatures on Grassmannians: moments, dimension reduction, and related
  topics}.
In \bbooktitle{Compressed Sensing and its Applications}
\bpages{235--259}.
\bpublisher{Springer}.
\end{bincollection}
\endbibitem

\bibitem{chafai2021}
\begin{barticle}[author]
\bauthor{\bsnm{Chafaï},~\bfnm{D}\binits{D.}}
(\byear{2021}).
\btitle{Random projections, marginals, and moments}.
\bjournal{URL https://djalil. chafai. net/docs/projections. pdf. Accessed on
  04.04}.
\end{barticle}
\endbibitem

\bibitem{cherian1941}
\begin{barticle}[author]
\bauthor{\bsnm{Cherian},~\bfnm{KC}\binits{K.}}
(\byear{1941}).
\btitle{A Bivariate Correlated Gamma-Type Distribution Function}.
\bjournal{Journal of the Indian Mathematical Society}
\bvolume{5}
\bpages{133--144}.
\end{barticle}
\endbibitem

\bibitem{chizat2021}
\begin{barticle}[author]
\bauthor{\bsnm{Chizat},~\bfnm{Lenaic}\binits{L.}}
(\byear{2021}).
\btitle{Sparse optimization on measures with over-parameterized gradient
  descent}.
\bjournal{Mathematical Programming}
\bpages{1--46}.
\end{barticle}
\endbibitem

\bibitem{chizat2018}
\begin{barticle}[author]
\bauthor{\bsnm{Chizat},~\bfnm{Lenaic}\binits{L.}} \AND
  \bauthor{\bsnm{Bach},~\bfnm{Francis}\binits{F.}}
(\byear{2018}).
\btitle{On the global convergence of gradient descent for over-parameterized
  models using optimal transport}.
\bjournal{arXiv preprint arXiv:1805.09545}.
\end{barticle}
\endbibitem

\bibitem{comon2010}
\begin{bbook}[author]
\bauthor{\bsnm{Comon},~\bfnm{Pierre}\binits{P.}} \AND
  \bauthor{\bsnm{Jutten},~\bfnm{Christian}\binits{C.}}
(\byear{2010}).
\btitle{Handbook of Blind Source Separation: Independent component analysis and
  applications}.
\bpublisher{Academic press}.
\end{bbook}
\endbibitem

\bibitem{comte2015}
\begin{barticle}[author]
\bauthor{\bsnm{Comte},~\bfnm{Fabienne}\binits{F.}} \AND
  \bauthor{\bsnm{{Genon-Catalot}},~\bfnm{Valentine}\binits{V.}}
(\byear{2015}).
\btitle{Adaptive {{Laguerre}} Density Estimation for Mixed {{Poisson}} Models}.
\bjournal{Electronic Journal of Statistics}
\bvolume{9}
\bpages{1113--1149}.
\end{barticle}
\endbibitem

\bibitem{constantine1996}
\begin{barticle}[author]
\bauthor{\bsnm{Constantine},~\bfnm{G.~M.}\binits{G.~M.}} \AND
  \bauthor{\bsnm{Savits},~\bfnm{T.~H.}\binits{T.~H.}}
(\byear{1996}).
\btitle{A {{Multivariate Faa}} Di {{Bruno Formula}} with {{Applications}}}.
\bjournal{Transactions of the American Mathematical Society}
\bvolume{348}
\bpages{503--520}.
\end{barticle}
\endbibitem

\bibitem{cramer1936}
\begin{barticle}[author]
\bauthor{\bsnm{Cram{\'e}r},~\bfnm{Harald}\binits{H.}} \AND
  \bauthor{\bsnm{Wold},~\bfnm{Herman}\binits{H.}}
(\byear{1936}).
\btitle{Some theorems on distribution functions}.
\bjournal{Journal of the London Mathematical Society}
\bvolume{1}
\bpages{290--294}.
\end{barticle}
\endbibitem

\bibitem{cvijovic2011}
\begin{barticle}[author]
\bauthor{\bsnm{Cvijovi{\'c}},~\bfnm{Djurdje}\binits{D.}}
(\byear{2011}).
\btitle{New Identities for the Partial {{Bell}} Polynomials}.
\bjournal{Applied Mathematics Letters}
\bvolume{24}
\bpages{1544--1547}.
\end{barticle}
\endbibitem

\bibitem{decastro2021}
\begin{barticle}[author]
\bauthor{\bsnm{{de Castro}},~\bfnm{Yohann}\binits{Y.}},
  \bauthor{\bsnm{Gadat},~\bfnm{S{\'e}bastien}\binits{S.}},
  \bauthor{\bsnm{Marteau},~\bfnm{Cl{\'e}ment}\binits{C.}} \AND
  \bauthor{\bsnm{Maugis},~\bfnm{Cathy}\binits{C.}}
(\byear{2021}).
\btitle{{{SuperMix}}: {{Sparse Regularization}} for {{Mixtures}}}.
\bjournal{The annals of Statistics}
\bvolume{49}
\bpages{1779--1809}.
\end{barticle}
\endbibitem

\bibitem{dinardo2009}
\begin{barticle}[author]
\bauthor{\bsnm{Di~Nardo},~\bfnm{E}\binits{E.}},
  \bauthor{\bsnm{Guarino},~\bfnm{G}\binits{G.}} \AND
  \bauthor{\bsnm{Senato},~\bfnm{D}\binits{D.}}
(\byear{2009}).
\btitle{A new method for fast computing unbiased estimators of cumulants}.
\bjournal{Statistics and Computing}
\bvolume{19}
\bpages{155--165}.
\end{barticle}
\endbibitem

\bibitem{nardo2011}
\begin{barticle}[author]
\bauthor{\bsnm{Di~Nardo},~\bfnm{Elvira}\binits{E.}},
  \bauthor{\bsnm{Guarino},~\bfnm{Giuseppe}\binits{G.}} \AND
  \bauthor{\bsnm{Senato},~\bfnm{Domenico}\binits{D.}}
(\byear{2011}).
\btitle{A new algorithm for computing the multivariate Faa di Bruno’s
  formula}.
\bjournal{Applied mathematics and computation}
\bvolume{217}
\bpages{6286--6295}.
\end{barticle}
\endbibitem

\bibitem{dozat2016}
\begin{barticle}[author]
\bauthor{\bsnm{Dozat},~\bfnm{Timothy}\binits{T.}}
(\byear{2016}).
\btitle{Incorporating nesterov momentum into adam}.
\end{barticle}
\endbibitem

\bibitem{dussap2021}
\begin{barticle}[author]
\bauthor{\bsnm{Dussap},~\bfnm{Florian}\binits{F.}}
(\byear{2021}).
\btitle{Anisotropic multivariate deconvolution using projection on the Laguerre
  basis}.
\bjournal{Journal of Statistical Planning and Inference}
\bvolume{215}
\bpages{23--46}.
\end{barticle}
\endbibitem

\bibitem{ferriero2011}
\begin{barticle}[author]
\bauthor{\bsnm{Ferriero},~\bfnm{Alessandro}\binits{A.}}
(\byear{2011}).
\btitle{A New Method for Modeling Dependence via Extended Common Shock Type
  Model}.
\bpages{20}.
\end{barticle}
\endbibitem

\bibitem{furman2017}
\begin{barticle}[author]
\bauthor{\bsnm{Furman},~\bfnm{Edward}\binits{E.}},
  \bauthor{\bsnm{Hackmann},~\bfnm{Daniel}\binits{D.}} \AND
  \bauthor{\bsnm{Kuznetsov},~\bfnm{Alexey}\binits{A.}}
(\byear{2017}).
\btitle{On {{Log}}-{{Normal Convolutions}}: {{An Analytical}}-{{Numerical
  Method With Applications}} to {{Economic Capital Determination}}}.
\bjournal{SSRN Electronic Journal}.
\end{barticle}
\endbibitem

\bibitem{heppes1956}
\begin{barticle}[author]
\bauthor{\bsnm{Heppes},~\bfnm{A}\binits{A.}}
(\byear{1956}).
\btitle{On the determination of probability distributions of more dimensions by
  their projections}.
\bjournal{Acta Mathematica Academiae Scientiarum Hungarica}
\bvolume{7}
\bpages{403--410}.
\end{barticle}
\endbibitem

\bibitem{james2008}
\begin{barticle}[author]
\bauthor{\bsnm{James},~\bfnm{Lancelot~F.}\binits{L.~F.}},
  \bauthor{\bsnm{Roynette},~\bfnm{Bernard}\binits{B.}} \AND
  \bauthor{\bsnm{Yor},~\bfnm{Marc}\binits{M.}}
(\byear{2008}).
\btitle{Generalized {{Gamma Convolutions}}, {{Dirichlet}} Means, {{Thorin}}
  Measures, with Explicit Examples}.
\bjournal{Probability Surveys}
\bvolume{5}
\bpages{346--415}.
\end{barticle}
\endbibitem

\bibitem{kim2019}
\begin{barticle}[author]
\bauthor{\bsnm{Kim},~\bfnm{Taekyun}\binits{T.}},
  \bauthor{\bsnm{Kim},~\bfnm{Dae}\binits{D.}} \AND
  \bauthor{\bsnm{Jang},~\bfnm{Gwan-Woo}\binits{G.-W.}}
(\byear{2019}).
\btitle{On {{Central Complete}} and {{Incomplete Bell Polynomials I}}}.
\bjournal{Symmetry}
\bvolume{11}
\bpages{288}.
\end{barticle}
\endbibitem

\bibitem{kingma2014}
\begin{barticle}[author]
\bauthor{\bsnm{Kingma},~\bfnm{Diederik~P}\binits{D.~P.}} \AND
  \bauthor{\bsnm{Ba},~\bfnm{Jimmy}\binits{J.}}
(\byear{2014}).
\btitle{Adam: A method for stochastic optimization}.
\bjournal{arXiv preprint arXiv:1412.6980}.
\end{barticle}
\endbibitem

\bibitem{klugman1999}
\begin{barticle}[author]
\bauthor{\bsnm{Klugman},~\bfnm{Stuart~A}\binits{S.~A.}} \AND
  \bauthor{\bsnm{Parsa},~\bfnm{Rahul}\binits{R.}}
(\byear{1999}).
\btitle{Fitting bivariate loss distributions with copulas}.
\bjournal{Insurance: mathematics and economics}
\bvolume{24}
\bpages{139--148}.
\end{barticle}
\endbibitem

\bibitem{laverny2021a}
\begin{barticle}[author]
\bauthor{\bsnm{Laverny},~\bfnm{Oskar}\binits{O.}},
  \bauthor{\bsnm{Masiello},~\bfnm{Esterina}\binits{E.}},
  \bauthor{\bsnm{Maume-Deschamps},~\bfnm{Véronique}\binits{V.}} \AND
  \bauthor{\bsnm{Rullière},~\bfnm{Didier}\binits{D.}}
(\byear{2021}).
\btitle{{Estimation of multivariate generalized gamma convolutions through
  Laguerre expansions.}}
\bjournal{Electronic Journal of Statistics}
\bvolume{15}
\bpages{5158 -- 5202}.
\bdoi{10.1214/21-EJS1918}
\end{barticle}
\endbibitem

\bibitem{mabon2017}
\begin{barticle}[author]
\bauthor{\bsnm{Mabon},~\bfnm{Gwenna{\"e}lle}\binits{G.}}
(\byear{2017}).
\btitle{Adaptive {{Deconvolution}} on the {{Non}}-Negative {{Real Line}}:
  {{Adaptive}} Deconvolution on {{R}}+}.
\bjournal{Scandinavian Journal of Statistics}
\bvolume{44}
\bpages{707--740}.
\end{barticle}
\endbibitem

\bibitem{mathai1982}
\begin{barticle}[author]
\bauthor{\bsnm{Mathai},~\bfnm{A.~M.}\binits{A.~M.}}
(\byear{1982}).
\btitle{Storage Capacity of a Dam with Gamma Type Inputs}.
\bjournal{Annals of the Institute of Statistical Mathematics}
\bvolume{34}
\bpages{591--597}.
\end{barticle}
\endbibitem

\bibitem{melman2018}
\begin{barticle}[author]
\bauthor{\bsnm{Melman},~\bfnm{Vadim~S}\binits{V.~S.}},
  \bauthor{\bsnm{Shablya},~\bfnm{Yuriy~V}\binits{Y.~V.}},
  \bauthor{\bsnm{Kruchinin},~\bfnm{Dmitry~V}\binits{D.~V.}} \AND
  \bauthor{\bsnm{Shelupanov},~\bfnm{Alexander~A}\binits{A.~A.}}
(\byear{2018}).
\btitle{Realization of a Method for Calculating Bell Polynomials Based on
  Compositae of Generating Functions}.
\bjournal{Journal of Informatics and Mathematical Sciences}
\bvolume{10}
\bpages{659--672}.
\end{barticle}
\endbibitem

\bibitem{miatto2019}
\begin{barticle}[author]
\bauthor{\bsnm{Miatto},~\bfnm{Filippo~M.}\binits{F.~M.}}
(\byear{2019}).
\btitle{Recursive Multivariate Derivatives of $e^{f(X_1,\ldots, X_n)}$ of
  Arbitrary Order}.
\bjournal{arXiv:1911.11722 [cs, math]}.
\end{barticle}
\endbibitem

\bibitem{miles2019}
\begin{barticle}[author]
\bauthor{\bsnm{Miles},~\bfnm{Justin}\binits{J.}},
  \bauthor{\bsnm{Furman},~\bfnm{Edward}\binits{E.}} \AND
  \bauthor{\bsnm{Kuznetsov},~\bfnm{Alexey}\binits{A.}}
(\byear{2021}).
\btitle{Risk Aggregation: {{A}} General Approach via the Class of {{Generalized
  Gamma Convolutions}}}.
\bjournal{Variance}
\bvolume{13:2}
\bpages{233--249}.
\end{barticle}
\endbibitem

\bibitem{mishkov2000}
\begin{barticle}[author]
\bauthor{\bsnm{Mishkov},~\bfnm{Rumen~L.}\binits{R.~L.}}
(\byear{2000}).
\btitle{Generalization of the Formula of {{Faa}} Di {{Bruno}} for a Composite
  Function with a Vector Argument}.
\bjournal{International Journal of Mathematics and Mathematical Sciences}
\bvolume{24}
\bpages{481--491}.
\end{barticle}
\endbibitem

\bibitem{moschopoulos1985}
\begin{barticle}[author]
\bauthor{\bsnm{Moschopoulos},~\bfnm{P.~G.}\binits{P.~G.}}
(\byear{1985}).
\btitle{The Distribution of the Sum of Independent Gamma Random Variables}.
\bjournal{Annals of the Institute of Statistical Mathematics}
\bvolume{37}
\bpages{541--544}.
\end{barticle}
\endbibitem

\bibitem{perez-abreu2012}
\begin{barticle}[author]
\bauthor{\bsnm{{P{\'e}rez-Abreu}},~\bfnm{Victor}\binits{V.}} \AND
  \bauthor{\bsnm{Stelzer},~\bfnm{Robert}\binits{R.}}
(\byear{2012}).
\btitle{A {{Class}} of {{Infinitely Divisible Multivariate}} and {{Matrix Gamma
  Distributions}} and {{Cone}}-Valued {{Generalised Gamma Convolutions}}}.
\bjournal{arXiv:1201.1461 [math, stat]}.
\end{barticle}
\endbibitem

\bibitem{perez-abreu2014}
\begin{barticle}[author]
\bauthor{\bsnm{{P{\'e}rez-Abreu}},~\bfnm{Victor}\binits{V.}} \AND
  \bauthor{\bsnm{Stelzer},~\bfnm{Robert}\binits{R.}}
(\byear{2014}).
\btitle{Infinitely Divisible Multivariate and Matrix {{Gamma}} Distributions}.
\bjournal{Journal of Multivariate Analysis}
\bvolume{130}
\bpages{155--175}.
\end{barticle}
\endbibitem

\bibitem{renyi1952}
\begin{barticle}[author]
\bauthor{\bsnm{R{\'e}nyi},~\bfnm{Alfr{\'e}d}\binits{A.}}
(\byear{1952}).
\btitle{On projections of probability distributions}.
\bjournal{Acta Math. Acad. Sci. Hungar}
\bvolume{3}
\bpages{131--142}.
\end{barticle}
\endbibitem

\bibitem{roynette2009}
\begin{barticle}[author]
\bauthor{\bsnm{Roynette},~\bfnm{Bernard}\binits{B.}},
  \bauthor{\bsnm{Vallois},~\bfnm{Pierre}\binits{P.}} \AND
  \bauthor{\bsnm{Yor},~\bfnm{Marc}\binits{M.}}
(\byear{2009}).
\btitle{A Family of Generalized Gamma Convoluted Variables}.
\end{barticle}
\endbibitem

\bibitem{smith2020}
\begin{barticle}[author]
\bauthor{\bsnm{Smith},~\bfnm{Kevin~D.}\binits{K.~D.}}
(\byear{2020}).
\btitle{A {{Tutorial}} on {{Multivariate}} $k$-{{Statistics}} and Their
  {{Computation}}}.
\bjournal{arXiv:2005.08373 [math, stat]}.
\end{barticle}
\endbibitem

\bibitem{smith1995}
\begin{barticle}[author]
\bauthor{\bsnm{Smith},~\bfnm{Peter~J}\binits{P.~J.}}
(\byear{1995}).
\btitle{A recursive formulation of the old problem of obtaining moments from
  cumulants and vice versa}.
\bjournal{The American Statistician}
\bvolume{49}
\bpages{217--218}.
\end{barticle}
\endbibitem

\bibitem{thorin1977}
\begin{barticle}[author]
\bauthor{\bsnm{Thorin},~\bfnm{Olof}\binits{O.}}
(\byear{1977}).
\btitle{On the Infinite Divisibility of the {{Pareto}} Distribution}.
\bjournal{Scandinavian Actuarial Journal}
\bvolume{1977}
\bpages{31--40}.
\end{barticle}
\endbibitem

\bibitem{thorin1977a}
\begin{barticle}[author]
\bauthor{\bsnm{Thorin},~\bfnm{Olof}\binits{O.}}
(\byear{1977}).
\btitle{On the Infinite Divisibility of the Lognormal Distribution}.
\bjournal{Scandinavian Actuarial Journal}
\bvolume{1977}
\bpages{121--148}.
\end{barticle}
\endbibitem

\bibitem{wendland2004}
\begin{bbook}[author]
\bauthor{\bsnm{Wendland},~\bfnm{Holger}\binits{H.}}
(\byear{2004}).
\btitle{Scattered data approximation}
\bvolume{17}.
\bpublisher{Cambridge university press}.
\end{bbook}
\endbibitem

\end{thebibliography}


\end{document}